\documentclass[final]{elsarticle}

\usepackage[margin=1.25in]{geometry}
\usepackage{amssymb}
\usepackage{amsthm,amsmath,amssymb,braket,amsfonts}
\usepackage{xspace}
\usepackage{float}
\usepackage{etoolbox}
\usepackage{algorithm}
\usepackage{algpseudocode}
\usepackage{cuted}
\usepackage{hyperref}
\usepackage[capitalise,nameinlink]{cleveref}
\crefname{equation}{}{}
\crefname{appendix}{}{}

\renewcommand{\vec}[1]{{\mathchoice
                {\mbox{\boldmath$\displaystyle{#1}$}}
                {\mbox{\boldmath$\textstyle{#1}$}}
                {\mbox{\boldmath$\scriptstyle{#1}$}}
                {\mbox{\boldmath$\scriptscriptstyle{#1}$}}}}
\newcommand{\mat}[1]{\mathbf{{#1}}}

\newcommand{\R}{\mathbb{R}}
\newcommand{\E}{\mathbb{E}}
\newcommand{\var}{\mathrm{var}}
\newcommand{\f}{\hat{f}}

\newcommand{\MATLAB}{\textsc{Matlab}}

\usepackage{lipsum}
\usepackage{subfigure}
\usepackage{graphicx}
\usepackage{epstopdf}

\ifpdf
  \DeclareGraphicsExtensions{.eps,.pdf,.png,.jpg}
\else
  \DeclareGraphicsExtensions{.eps}
\fi

\usepackage{enumitem}
\setlist[enumerate]{leftmargin=.5in}
\setlist[itemize]{leftmargin=.5in}

\newtheorem{theorem}{Theorem}[section]

\newtheorem{lemma}[theorem]{Lemma}
\newtheorem{proposition}[theorem]{Proposition}

\journal{arXiv}

\begin{document}

\begin{frontmatter}



\title{Extreme learning machines for variance-based global sensitivity analysis\tnoteref{thanks,funding}}
\tnotetext[thanks]{Submitted to the editors \today.}
\tnotetext[funding]{Supported in part by US National Science Foundation grants DMS \#1745654 and DMS \#1953271.}

\author[label1]{John Darges\corref{cor1}}\ead{jedarges@ncsu.edu}

\author[label1]{Alen Alexanderian}\ead{alexanderian@ncsu.edu}

\author[label1,label2]{Pierre Gremaud}\ead{gremaud@ncsu.edu}

\affiliation[label1]{organization={Department of Mathematics, North Carolina State University},
city={Raleigh},
state={NC},
 country={USA}}
 \affiliation[label2]{organization={The Graduate School, North Carolina State University},
city={Raleigh},
state={NC},
country={USA}}

\cortext[cor1]{Corresponding author}

\begin{abstract}
Variance-based global sensitivity analysis (GSA) can provide a wealth of
information when applied to complex models. A well-known Achilles' heel of this
approach is its computational cost which often renders it unfeasible in
practice. An appealing alternative is to analyze instead the sensitivity of a
surrogate model with the goal of lowering computational costs while maintaining
sufficient accuracy. Should a surrogate be ``simple" enough to  be amenable to
the analytical calculations of its Sobol' indices, the cost of GSA is
essentially reduced to the construction of the surrogate.  We propose a new
class of sparse weight Extreme Learning Machines (SW-ELMs) which, when
considered as surrogates in the context of GSA, admit analytical formulas for
their Sobol' indices and, unlike the standard ELMs, yield accurate
approximations of these indices. The effectiveness of this approach is
illustrated through both traditional benchmarks in the field and on a chemical
reaction network. 
\end{abstract}


\end{frontmatter}


\section{Introduction}\label{sec:intro}
 A key insight into the behavior of  a generic model of the form 
\begin{equation}\label{equ:model}
    y = f(\vec{x}), \quad \vec{x} \in \R^d, y \in \R,
\end{equation}
is to understand the impact of the uncertainty in the entries of $\vec{x} = (x_1, \dots, x_d)$
on the uncertainty in the model output $y$; this can be addressed by
performing global sensitivity analysis
(GSA)~\cite{Saltelli08,IoossLeMaitre15,IoossSaltelli17}.  We focus on variance
based GSA~\cite{SalSob95,Sobol01,PrieurTarantola17} whereby one seeks to quantify the relative contributions of the entries of $\vec{x}$ to the variance of $f$. Specifically, we rely on the Sobol' indices~\cite{Sobol01} and assume that  the entries of $\vec{x}$ in~\cref{equ:model} can be regarded as
independent uniformly distributed random variables\textemdash which is common in
applications.

Computing Sobol' indices generally  involves a costly sampling procedure and
multiple applications of Monte Carlo integration. This is not feasible if the
input dimension $d$ is large, a common situation in engineering applications. 
Often, an initial parameter screening procedure can be applied to reduce the
input dimension, before a more detailed variance-based GSA is conducted; see
e.g.,~\cite{HartGremaudDavid19}.  Also, specific problem structures may enable
input dimension reduction in creative ways. An interesting example is the
work~\cite{Tabandeh22}, where the authors consider hierarchical models and
devise a multi-level method for dimensionality reduction coupled with variable
grouping to reduce the number of function evaluations.  Nonetheless, when
evaluating the model $f$ is prohibitively expensive, Monte Carlo methods for
computing Sobol' indices are not viable.  In such cases, a widely used approach
is to construct a surrogate model $\hat{f} \approx f$ whose Sobol' indices can
be computed efficiently~\cite{GratietMarelliSudret17,Sargsyan2017}.  

A number of surrogate models have been proposed for accelerating variance based 
GSA. These include polynomial chaos expansions (PCE)~\cite{Sudret08,Crestaux09},
multivariate adaptive regression splines
(MARS)~\cite{friedman91,HartAlexanderianGremaud17}, Gaussian
processes~\cite{Marrel2IoossLaurentEtAl09,OakleyOHagan04,JinChenSudjianto04},
Bayesian adaptive regression trees (BART)~\cite{Horiguchi21}, random
forests~\cite{Antoniadis21}, support vector machines~\cite{Steiner19},
artificial neural networks (ANN)~\cite{Fock14}, autoregressive
models~\cite{Datteo18}, and radial basis functions (RBF)~\cite{Leak14}.  
See~\cite{GratietMarelliSudret17,Kai20} for an overview of surrogate-based methods for GSA. 

In surrogate-based approaches, the user builds the surrogate model and estimates
the surrogate's Sobol' indices by Monte Carlo integration. This process
introduces two sources of error in the Sobol' indices. The first one is due to
the approximation error of the surrogate and the second one due to the error in
estimating the surrogate's Sobol' indices. Sampling of the surrogate is
often considerably cheaper than sampling the underlying model. However, this
sampling can still be costly.  In the best case, analytic formulas for the
surrogate's Sobol' indices are known. Then, the issues caused by Monte Carlo
sampling disappear altogether. Certain statistical surrogates such as Gaussian
processes~\cite{OakleyOHagan04,JinChenSudjianto04} and BART~\cite{Horiguchi21},
offer analytic Sobol' indices along with accompanying uncertainty bounds.  The
statistical nature of these surrogates, however, means sampling is needed to
compute Sobol' indices. This could be potentially costly,
especially compared to certain numerical surrogate models for which Sobol'
indices can be computed analytically for practically no additional
cost~\cite{GratietMarelliSudret17}.  
Examples of such approaches include PCEs~\cite{Sudret08,Crestaux09} or
RBFs~\cite{Wu16,Wu19}. 

In the case of PCE surrogates, analytic formulas for Sobol' indices derived from
the PC coefficients with no additional sampling
needed~\cite{Sudret08,Crestaux09}. Unsurprisingly, PCE has become a workhorse
method for computing Sobol' indices and uncertainty analysis in general.  This
has spawned a number of methods for computing and using PCEs for uncertainty
analysis. Non-intrusive project-based approaches involve computing the PC
coefficients using quadrature.  For high-dimensional inputs, this is untenable
and motivates regression-based approaches.  A major challenge with PCE is that
the PC basis can grow to gargantuan sizes as the degree of the expansion is
increased.  This necessitates sparse regression methods to reduce the size of
basis~\cite{Blatman10}. A survey of these methods can be found in~\cite{Luthen21}.  We also
point to~\cite{AlexanderianGremaudSmith20,Ehre20,Zhou20,Luthen22,Almohammadi22}, for recent developments regarding
PCE-based approaches for uncertainty quantification.  Overall, PCEs provide an
efficient computational tool for fast uncertainty analysis in problems with low
to moderate parameter dimensions.

In this article, we present a new approach which uses extreme learning machine
(ELM) surrogates---a class of neural networks---for which Sobol' indices can be
computed analytically. Our contribution expands the limited toolbox of fast
surrogate models with analytic formulas for Sobol' indices. Having ready access
to such analytic formulas enables the users to compute these indices without any
costly Monte Carlo sampling procedures.  The existing surrogate modeling
approaches that have this feature each have their strengths and weaknesses,
making them suitable for specific applications of interest. When studying black
box models, it is seldom obvious which surrogate is the best choice. In
such applications where one does not know the true values of the sensitivity
indices, it is beneficial to have access to different tools to compute the
indices. This can be useful for confirming the
reliability of one's results. An example of this is~\cite{Zhang22}, where
comparing the Sobol' indices of PCE and ANN surrogates provides verifiable GSA.
Notably, however, the ANN surrogate requires much more time to produce GSA
results.

ELMs are a class of single layer neural
networks~\cite{HuangELM,Huang11,Stewart21} where, unlike traditional neural
networks, one draws weights and biases of the hidden layer randomly.
Consequently, training an ELM amounts to solving a linear least squares problem
to estimate the output layer weights; see \cref{sec:ELMs} for more details.
Though the idea of randomly sampling weights may appear far-fetched, ELMs are a
theoretically sound approach. As a trade-off, ELMs may require more hidden layer
neurons than traditional neural networks to achieve the same results.  
Generally, ELMs have proven to be a powerful tool in a wide range of 
applications~\cite{Huang15,Wang22}.

Our motivation for considering ELMs for GSA  is twofold
\begin{itemize}
\item  they are simple and inexpensive to train: only a linear least squares problem needs to be solved,
\item through a judicious choice of the activation function, analytic formulas for
the Sobol' indices can be derived;  see~\cref{sec:gsa}.
\end{itemize}
Also, ELMs provide a flexible surrogate modeling
framework: they do not make strong demands on the distribution of the inputs or
the regularity of the input-output map. 
Additionally, ELMs are well suited to tackle nonlinear and high-dimensional models. 
The use of neural networks for surrogate-based uncertainty quantification has
been the object of recent research efforts. Among many others,
\cite{Nagawkar21,Walzberg21,Li19,Kapusuzoglu21,Ye21,Zhao23}  have for instance
used neural network surrogates for GSA. To our knowledge, however, the present
work is the first of its kind to derive analytic formulas for Sobol' indices of
neural network surrogates. 

%
% sparsity control
%
To be used reliably for variance based GSA, the surrogate $\hat f$ 
must capture some key structural properties of the exact model $f$. 
To illustrate this point, consider~\cref{equ:model} with $\vec{x}\in\R^3$ and assume $f$ has mean zero.
The ANOVA decomposition of $f$ is given by 
\begin{equation}\label{equ:anova3d}
    f(x_1,x_2,x_3) = \sum_{i=1}^3 f_i(x_i) + f_{12}(x_1,x_2) + f_{13}(x_1, x_3) 
   + f_{23}(x_2,x_3) + f_{123}(x_1,x_2,x_3),
\end{equation}
where $f_i = \mathbb{E}(f | x_i)$, $f_{ij} = \mathbb{E}(f | x_i,x_j)  - f_i - f_j$, and 
$f_{123} = f - \sum_{i=1}^3 f_i - \sum_{i < j} f_{ij}$. A key observation is that the 
variance of $f$ can be decomposed as the sum of the variances of the individual terms in the ANOVA.
This makes it possible to  quantify the contributions of the individual inputs
(or a group of inputs) to the total variance of $f$ and leads to definition of
Sobol' indices.  To provide accurate estimates of Sobol' indices, a surrogate
$\hat{f}$ must emulate the main effect terms $\{f_i\}_{i=1}^d$ as well as the
higher order interaction terms in~\cref{equ:anova3d}.  As discussed
in~\cref{sec:sparsity}, standard ELMs may fail to correctly capture the impact
of variable interactions on output variance, leading to inaccurate Sobol'
indices. We resolve this shortcoming by introducing sparsity to the hidden layer
weights. The sparsification technique is most closely related to network
pruning, a common method that improves training efficiency for neural
networks~\cite{Engelbrecht01,Blalock2020}. This unique feature of our method
preserves the speed and simplicity of ELM while also allowing the surrogate to
better incorporate model behavior.  We demonstrate the efficiency of the
proposed approach on standard benchmark problems, an application from
biochemistry, and a high-dimensional problem; see~\cref{sec:numerics}.

\section{Extreme learning machines}\label{sec:ELMs}
A single layer neural network (SLNN) with $n$ hidden layer neurons and a scalar output has the form
\begin{equation}\label{equ:onelayerNN}
   g(\vec{x}) = 
\sum_{j=1}^n \beta_j 
   \phi(\vec{w}_j^\top \vec{x} + b_j), 
\end{equation}
where $\vec{w}_j \in \mathbb R^d$, $b_j \in \mathbb R$, $j = 1, \ldots, n$, are the weights and biases of the 
hidden layer, $\beta_j$, $j = 1, \ldots, n$, are the weights of the output layer, 
and $\phi$ is the hidden layer activation function. For a given number of neurons $n$ and a given activation function $\phi$, we set 
\begin{equation}
   \mathcal M_n(\phi) = \Big\{ \sum_{j=1}^n \beta_j \phi(\vec{w}_j^\top \vec{x} + b_j) 
   \, : \, b_j, \beta_j \in \mathbb R, \vec{w}_j \in \mathbb R^d\Big\}.
   \end{equation}
It is known that SLNNs are universal approximators since, when using non-polynomial continuous activation functions, they are dense in $\mathcal C(\mathbb R^d)$~\cite{Pinkus99}; more precisely, if $\phi \in \mathcal C (\mathbb R)$ is not polynomial and $f \in \mathcal C(K)$ where $K$ is a compact subset of $\mathbb R^d$ then, for any $\epsilon >0$, there exists $n$ and $g_n \in \mathcal M_n(\phi)$ such that 
\begin{equation}\label{equ:univap}
\max _{\vec{x} \in K} | f(\vec x) - g_n(\vec{x})| < \epsilon.
\end{equation}
In addition, SLNNs are universal approximators on $L^p(\R^d)$, for $p\in[1,\infty)$, when the activation function is not a polynomial (almost everywhere)~\cite{Leshno93,Huang11}.

The standard approach to train (\ref{equ:onelayerNN}) is  to  determine all hidden layer weights and biases and output weights by solving a nonlinear least squares problem. 

In an extreme learning machine (ELM)~\cite{HuangELM},  the weight vectors
$\vec{w}_j$ and biases $b_j$, $j = 1, \ldots, n$, of the hidden layer are not
determined as part of a regression process but rather are {\em chosen
randomly}. Training an ELM then only involves determining the output layer
weights $\{\beta_j\}_{j = 1}^n$ by  solving a {\em linear} least squares
problem.  Remarkably, even though it bypasses training the hidden layer
weights---and replaces a costly nonlinear least squares problem by a
linear one---ELMs are (almost surely) universal approximators on $L^2(\R^d)$.  Namely, 
by Theorem 2.4 in~\cite{Huang11}, ELMs have a universal
approximation property (\ref{equ:univap}) on $L^2(\R^d)$ provided 
(i) 
the hidden layer weights and biases are sampled from a continuous probability distribution; (ii) 
the activation function is piecewise continuous and not a polynomial (almost everywhere); and (iii)
the output weights are determined by ordinary least squares.

\subsection{Computing ELM surrogates}
Let  
\[
\begin{aligned}
\mat{W} &= [\begin{matrix}\vec{w}_1 & \vec{w}_2 & \cdots & \vec{w}_n\end{matrix}]^\top, \\
\vec{b} &= [\begin{matrix} b_1 & b_2 & \cdots & b_n\end{matrix}]^\top,\\ 
\vec{\beta} &= [\begin{matrix} \beta_1 & \beta_2 & \cdots & \beta_n\end{matrix}]^\top;
\end{aligned}
\] 
the ELM~\eqref{equ:onelayerNN} takes the form
\begin{equation}\label{equ:onelayermatrix}
   g(\vec{x}) = 
\vec{\beta}^\top \phi(\mat{W} \vec{x} + \vec{b}),
\end{equation}
where $\phi$ is understood to act componentwise. The entries of the weight matrix $\mat{W}$ and bias vector $\vec{b}$ 
are sampled independently from a continuous probability distribution $\mathcal{D}$. It is important to note
that this is done before the training, and a fixed realization of $\mat{W}$ and $\vec{b}$ is used in all subsequent
computations.

To construct a surrogate, the model $f$ is sampled at $\{\vec{x}_i\}_{i=1}^m$, yielding $y_i = f(\vec{x}_i)$, $i=1, \dots, m$. We then find $\vec{\beta}$ by solving
\begin{equation}\label{equ:lsq}
   \min_{\vec\beta \in \R^n } \| \mat{H}\vec\beta - \vec{y} \|_2,
\end{equation}
where $\vec{y} = [\begin{matrix} y_1 & \cdots & y_m\end{matrix}]^\top$  and
$
H_{ij} = \phi(\vec{w}_j^\top \vec{x}_i + b_j)$ for 
$i,j \in \{1,\ldots, m\} \times \{1, \ldots, n\}$.

In our computations, we implement a regularized least squares problem to control the magnitude of the solution $\vec\beta$: 
\begin{equation}\label{equ:rlsq}
    \min_{\vec\beta} \frac12 \| \mat{H} \vec\beta - \vec{y} \|_2^2 + \frac\alpha2 \| \vec\beta\|_2^2,
\end{equation}
with $\alpha >0$, the solution of which is
\[
    \vec\beta^* = (\mat{H}^\top \mat{H} + \alpha \mat{I})^{-1} \mat{H}^\top \vec{y}.
\]
The regularization parameter $\alpha$ can be selected using the L-curve method or generalized cross validation (GCV)~\cite{Hansen10}. We construct the hidden layer weight matrix and bias vector by sampling individual values from the standard normal distribution and use Latin hypercube sampling (LHS) to sample points in the training set $\{ \vec{x}_i \}_{i=1}^m$.

\section{Global sensitivity analysis using ELMs}\label{sec:gsa}
For each entry $x_k$, $k=1, \dots, d$, the first-order Sobol' index $S_k$ and total Sobol' index  $S_k^{\text{tot}}$\ of a surrogate $y = \hat f(\vec{x})$ are defined as
\begin{equation}\label{equ:sobol}
S_k = \frac{\text{var}(\hat{f}_k)}{\text{var}(\hat{f})}, \quad S_k^{\text{tot}} = 1 - \frac{\text{var}(\E (\hat{f}|x_l,\ l\neq k))}{\text{var}(\hat{f})},
\end{equation}
where $\hat f_k(x_k) = \mathbb E(\hat f|x_k) - \mathbb E(\hat f)$; see e.g., \cite{PrieurTarantola17}. 
Henceforth, we assume that the entries of the input vector $\vec{x}$ are independent and uniformly distributed
on the interval $[0, 1]$ and, therefore, the input domain is $[0, 1]^d$.
It is straightforward to 
extend the proposed approach to the case where entries of $\vec{x}$ are independent 
uniformly distributed random variables on arbitrary closed and bounded intervals.

To obtain analytic formulas for the Sobol' indices while avoiding Monte Carlo approximations, we must be able to easily compute the mean, variance, and partial variances of the surrogate. By carefully choosing the activation function, we can design $\hat f$ as an ELM which (i) suits the above requirement and (ii)  preserves the universal approximation property.  Activation functions traditionally used in machine learning are not suited for our purpose as the corresponding calculations in (\ref{equ:sobol}) become impractical or impossible.  Instead, we choose an exponential activation function 
\[
\phi(t)=e^t, 
\]
which results in 
\begin{equation}\label{equ:elm}
\hspace{-1.25mm}
\hat{f}(\vec x) = \sum_{j=1}^n\beta_j e^{\vec{w}_j^\top \vec{x} + b_j} = \sum_{j=1}^n\Big(\beta_je^{b_j}\prod_{l=1}^de^{w_{j,l}x_l}\Big).
\end{equation}
Since $\phi$ is a continuous non-polynomial function, the results cited in~\cref{sec:ELMs} apply and ELMs constructed this way still have the universal approximation property. 

We now compute  the first and second moments of (\ref{equ:elm}).
\begin{lemma}\label{lma1}
%Assume $f:[0,1]^d\to\R$; 
The mean and variance of the ELM~\eqref{equ:elm} are
\[
\E(\hat{f}) = \sum_{j=1}^n\Big(\beta_je^{b_j}\prod_{l=1}^d\epsilon(w_{j,l})\Big)
\]
and
\[
\var(\hat{f}) = 
\sum_{j,i=1}^n\beta_j\beta_ie^{b_j+b_i}\Big(\prod_{l=1}^d\epsilon(w_{j,l}+w_{i,l})-\prod_{r=1}^d\epsilon(w_{j,r})\epsilon(w_{i,r}\Big),
\]
where $\epsilon(t) = \left\{\begin{array}{ccc}
\frac{e^t-1}{t}, && t\neq 0 \\
1, && t = 0
\end{array}\right.$
\end{lemma}
\begin{proof}
See~\cref{sec:derivation}.
\end{proof}

Using~\cref{lma1}, we obtain analytic expressions for 
the Sobol' indices of \eqref{equ:elm}. 
\begin{proposition}\label{prop1}
The first-order and total Sobol' indices of the ELM~\eqref{equ:elm} are given 
by 
\begin{equation}\label{equ:reg}
S_k=\frac{1}{\var(\hat{f})}\sum_{j,i=1}^n\beta_j\beta_ie^{b_j+b_i}\big(\epsilon(w_{j,k}+w_{i,k})-\epsilon(w_{j,k})\epsilon(w_{i,k})\big)\prod_{l\neq k}\epsilon(w_{j,l})\epsilon(w_{i,l})
\end{equation}
and
\begin{equation}\label{equ:tot}
S_k^{\mathrm{tot}}=1- \frac{1}{\var(\hat{f})}\sum_{j,i=1}^n\beta_j\beta_ie^{b_j+b_i}\epsilon(w_{j,k})\epsilon(w_{i,k})\Big(\prod_{r\neq k}\epsilon(w_{j,r}+w_{i,r})-\prod_{l\neq k}\epsilon(w_{j,r})\epsilon(w_{i,l})\Big),
\end{equation}
respectively, for $k = 1, \ldots, d$.
\end{proposition}
\begin{proof}
See~\cref{sec:derivation}.
\end{proof}

\section{Sparse weight ELMs}\label{sec:sparsity}

ELM surrogates \eqref{equ:onelayermatrix} are constructed by  randomly sampling the hidden layer weight matrix $\mathbf{W}$ and bias vector $\vec{b}$  from the standard normal distribution. As pointed out in Section~\ref{sec:ELMs}, an ELM surrogate $\hat f$ can be found to satisfy 
\begin{equation}\label{equ:approxf}
\hat  f \approx f
\end{equation}
as accurately as desired. Our goal is however not so much to construct $\hat f$ satisfying \eqref{equ:approxf} within a given tolerance but rather to construct $\hat f$ such that 
\begin{equation}\label{equ:approxs}
S(\hat  f ) \approx S(f),
\end{equation}
where $S$ stands here for any of the Sobol' indices from~\cref{prop1} and $S(\hat f)$ can be computed at low cost. 
To ensure that the Sobol' indices computed from $\hat f$ yield a reliable description of the variable interactions,
we thus have to construct a surrogate $\hat f$ that not only satisfies \eqref{equ:approxf} but also 
\begin{equation}\label{equ:approxu}
\hat f_u \approx f_u, \quad \mbox{for any } u \subset \{1, \dots, d\},
\end{equation}
where $f_u$ and $\hat f_u$ are the terms corresponding to the subset $u$ in the ANOVA decomposition of $f$ and $\hat f$ respectively.  

\subsection{ELMs and variable interactions}\label{sec:inter}
Simple examples show that the amount of variable interactions in $f$ plays a key role in either achieving or failing to achieve \eqref{equ:approxs}. As an illustration, consider the following parametrized function 
\begin{equation}\label{equ:bnch_int}
   f_\delta(\vec{x})=\sum_{i=1}^dx_i+\delta\prod_{j=1}^d(1+x_j),\quad\vec{x}\in[0,1]^d,
\end{equation}
where $\delta$ controls the amount of variable interactions. When $\delta = 0$, $f_\delta$ is fully additive: there are no variable interactions  and thus the first-order and corresponding total Sobol indices are equal. For $\delta >0$, there are  interactions between the variables ithat increase in importance as $\delta$ increases. By construction, we also observe that the individual  entries $x_i$, $i=1, \dots, d$, all share the same Sobol' indices; this feature is not present when considering the ELM surrogate $\hat f$ instead of $f$ because of approximation errors.  
The first-order and total Sobol' indices of \eqref{equ:bnch_int} can be found  analytically and are given in~\cref{sec:sobolinds}. Therefore, by varying $\delta$, we can explore how  variable interactions affect accuracy in both \eqref{equ:approxf} and \eqref{equ:approxs}.

\begin{figure}[h!!]
    \centering
	\subfigure{\includegraphics[width=0.49\textwidth]{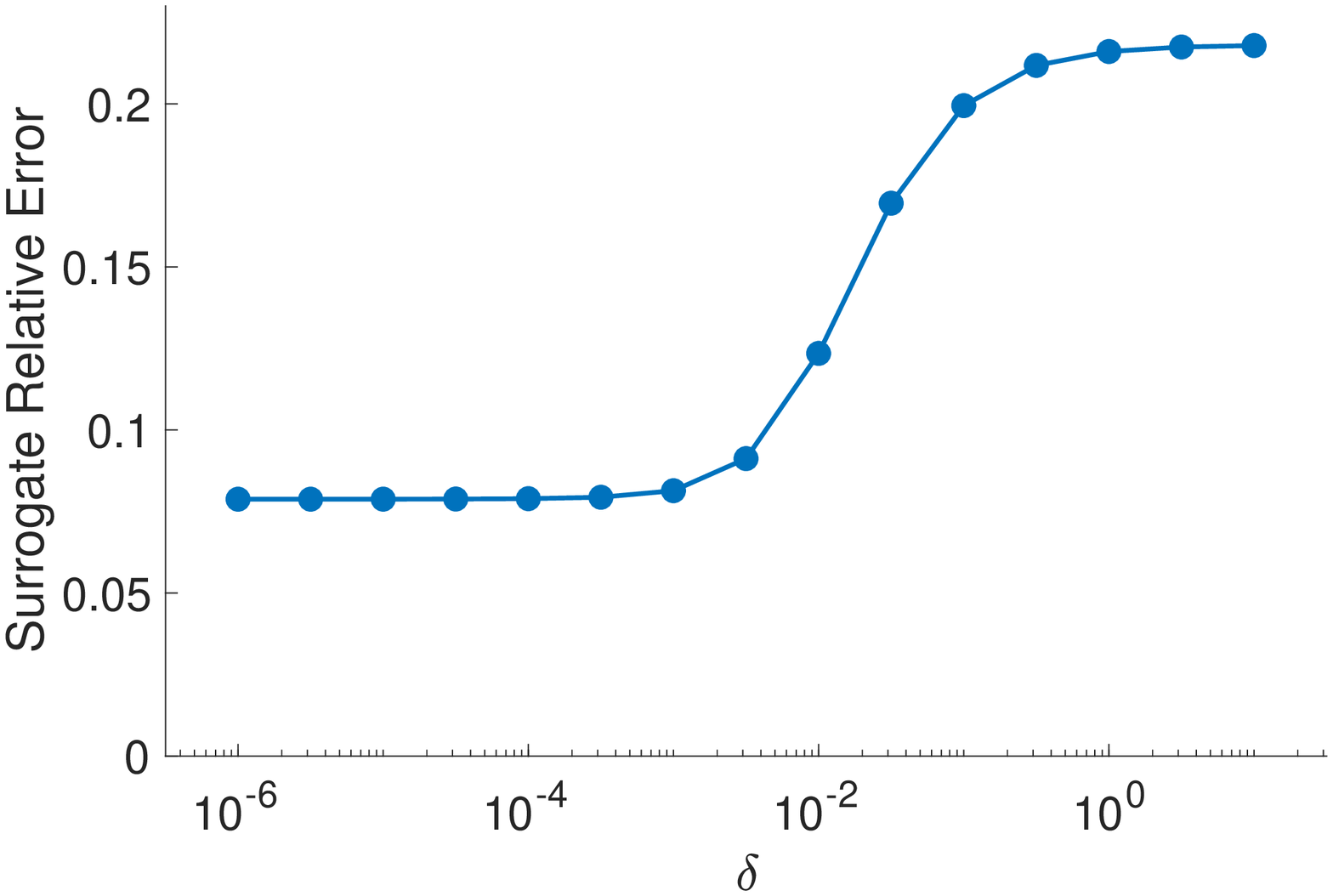}} 
   \subfigure{\includegraphics[width=0.49\textwidth]{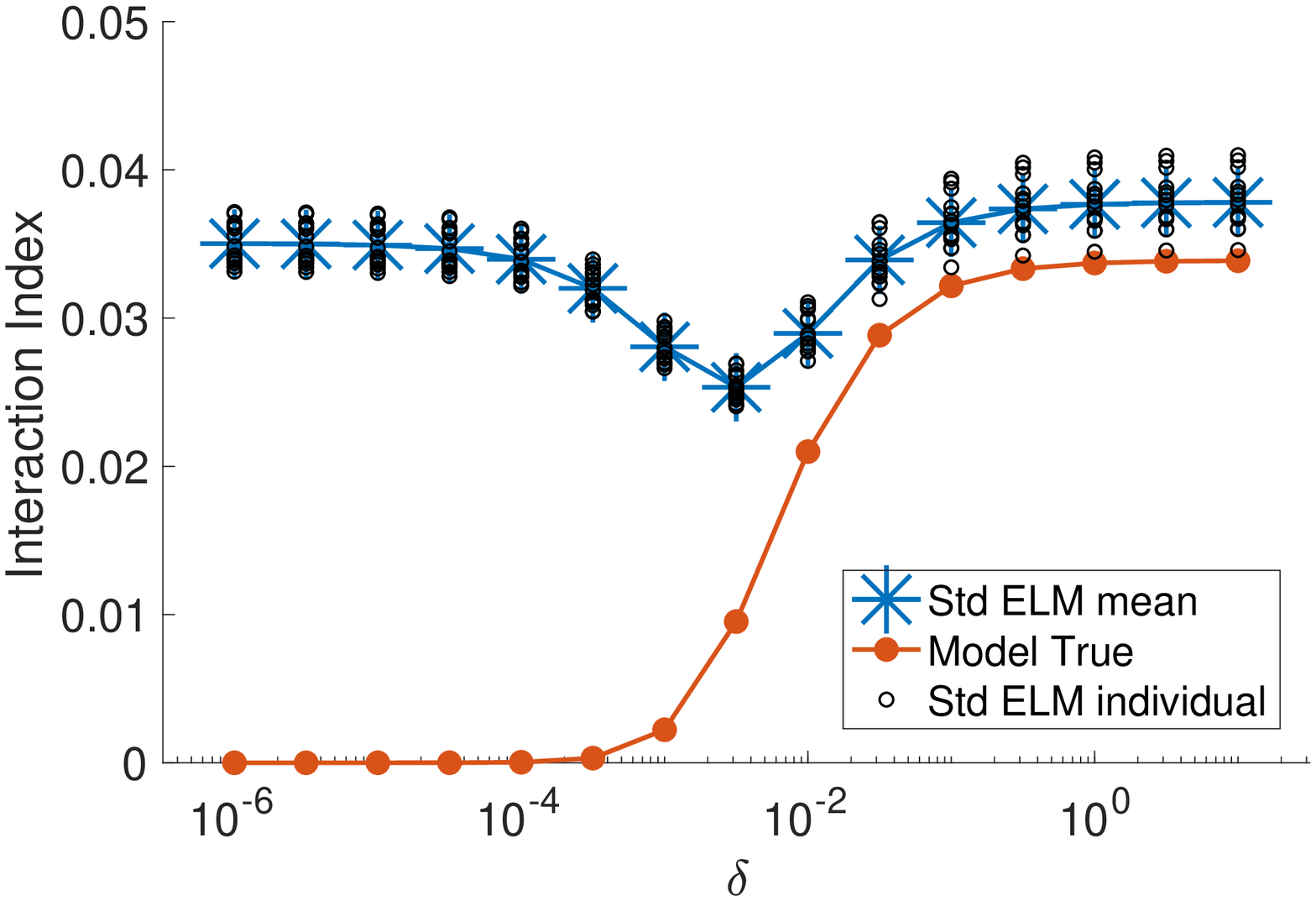}}
    \subfigure{\includegraphics[width=0.49\textwidth]{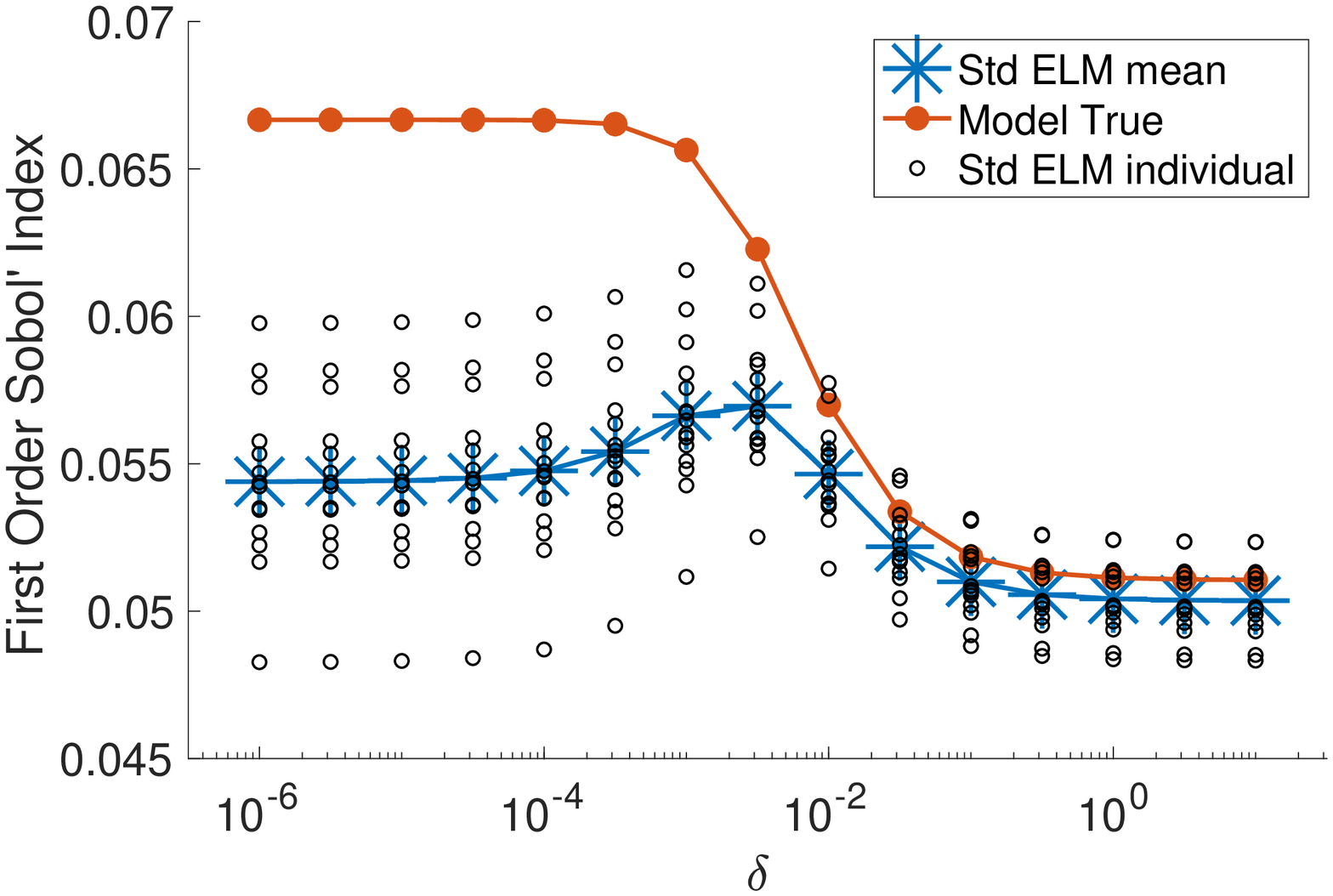}}
     \subfigure{\includegraphics[width=0.49\textwidth]{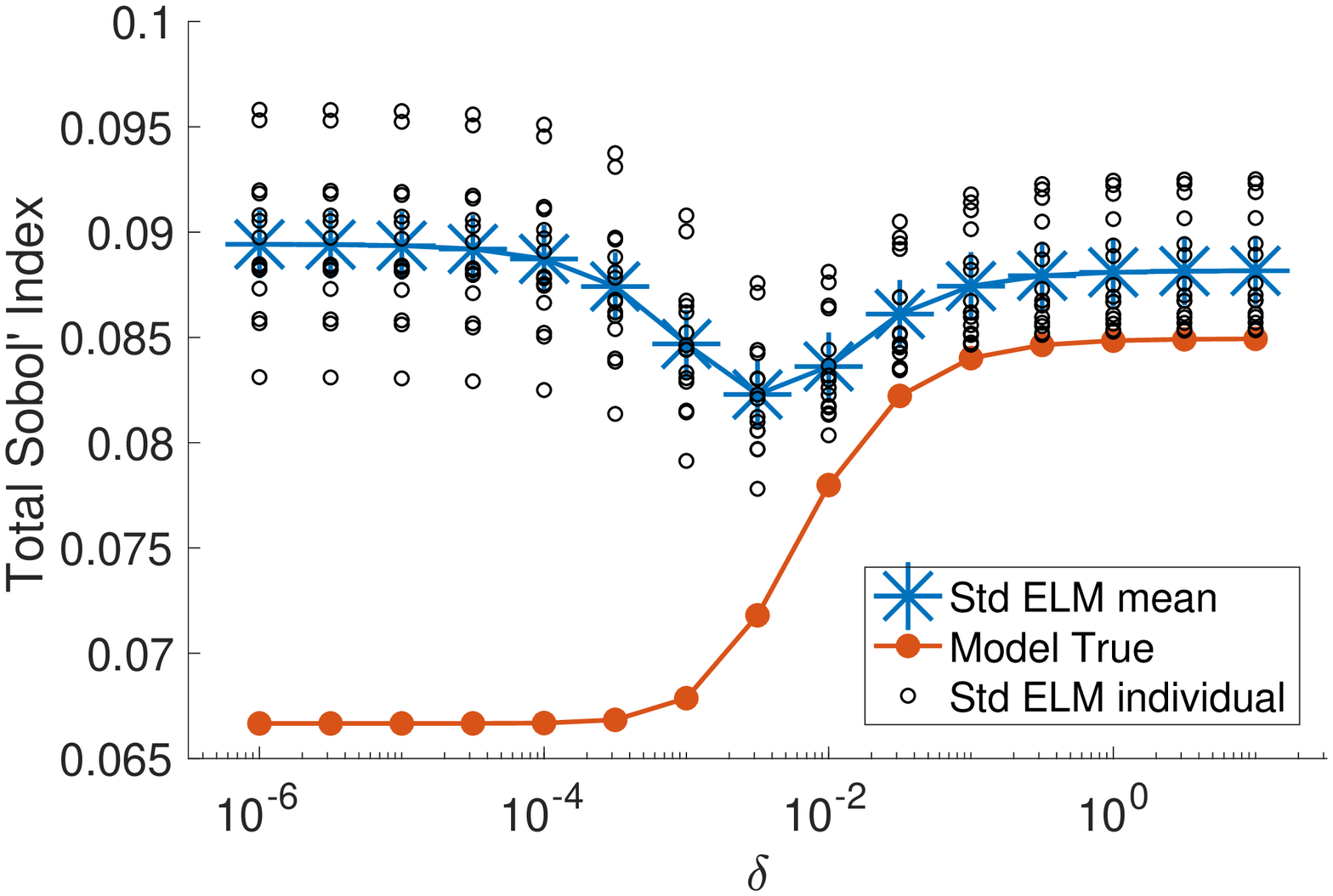}}
 \caption{Example~\eqref{equ:bnch_int} with $d=15$. Top left: Surrogate relative error; bottom left: first-order Sobol' indices $S_i(f)$ and $S_i(\hat f)$ of  $x_i$, $i=1,\dots,15$ and corresponding mean; bottom right: total Sobol' indices $S^\mathrm{tot}_i(f)$ and $S^\mathrm{tot}_i(\hat f)$  of $x_i$, $i=1,\dots,15$,  and corresponding mean; top right: interaction indices $I_i:=S^\mathrm{tot}_i-S_i$,  $i=1,\dots,15$, for both $f$ and $\hat f$.} 
\label{fig:bench_int}
\end{figure}

To approximate the Sobol' indices of \eqref{equ:bnch_int} with $d=15$, we construct ELMs of the form~\eqref{equ:elm}  as follows:

\begin{algorithmic}[1]
\State collect $m=900$ training points sampled by LHS
\State construct ELMs with $n=300$ neurons
\State compute the   surrogate relative error
\begin{equation}\label{equ:Esurr}
E_\mathrm{surr}=\frac{1}{\sqrt{\sum_{j=1}^sy_j^2}}\sqrt{\sum_{j=1}^s(\hat{f}(\vec{x}_j)-y_j)^2},\quad y_j=f(\vec{x}_j)\quad j=1,\dots,s,
\end{equation}
using $s=1000$ validation points sampled by LHS
\State select a regularization parameter $\alpha=10^{-3}$ by the L-curve method
\State compute Sobol' indices of the ELM surrogate using~\eqref{equ:reg} and~\eqref{equ:tot}
\end{algorithmic}

\cref{fig:bench_int} (top left) illustrates that, unsurprisingly, the surrogate relative error\textemdash which measures the discrepancy in \eqref{equ:approxf}\textemdash increases as the amount of interactions increases with $\delta$. Interestingly, the situation is reversed when looking at the accuracy of the Sobol' indices through \eqref{equ:approxs}. Indeed, \cref{fig:bench_int} (bottom left) shows that while $S_i(\hat f)$  is a close approximation of $S_i(f)$, for any $i=1,\dots, d$, at the larger values of $\delta$, the relative accuracy of the corresponding approximations decreases for smaller levels of interaction (i.e., smaller values of $\delta$). The variance of the $S_i(\hat f)$'s, which again should ideally all have the value $S_1(f)=\dots = S_d(f)$, also increases with smaller values of $\delta$.  The total Sobol' indices and their approximations largely behave in similar fashion, see  \cref{fig:bench_int} (bottom right). This example shows that ELM based surrogates have a tendency to overestimate variable interactions, a point made clear in \cref{fig:bench_int} (top right) where the interaction indices $I_i:=S^\mathrm{tot}_i-S_i$, $i=1, \dots, d$ are considered.

These results demonstrate that, for variance based GSA, standard ELM may fail to effectively adapt when applied to models featuring different degrees of contribution from interaction terms to the output variance.

\subsection{Hidden layer weight sparsification}

\begin{figure}[h!!]
    \centering	\subfigure{\includegraphics[width=0.49\textwidth]{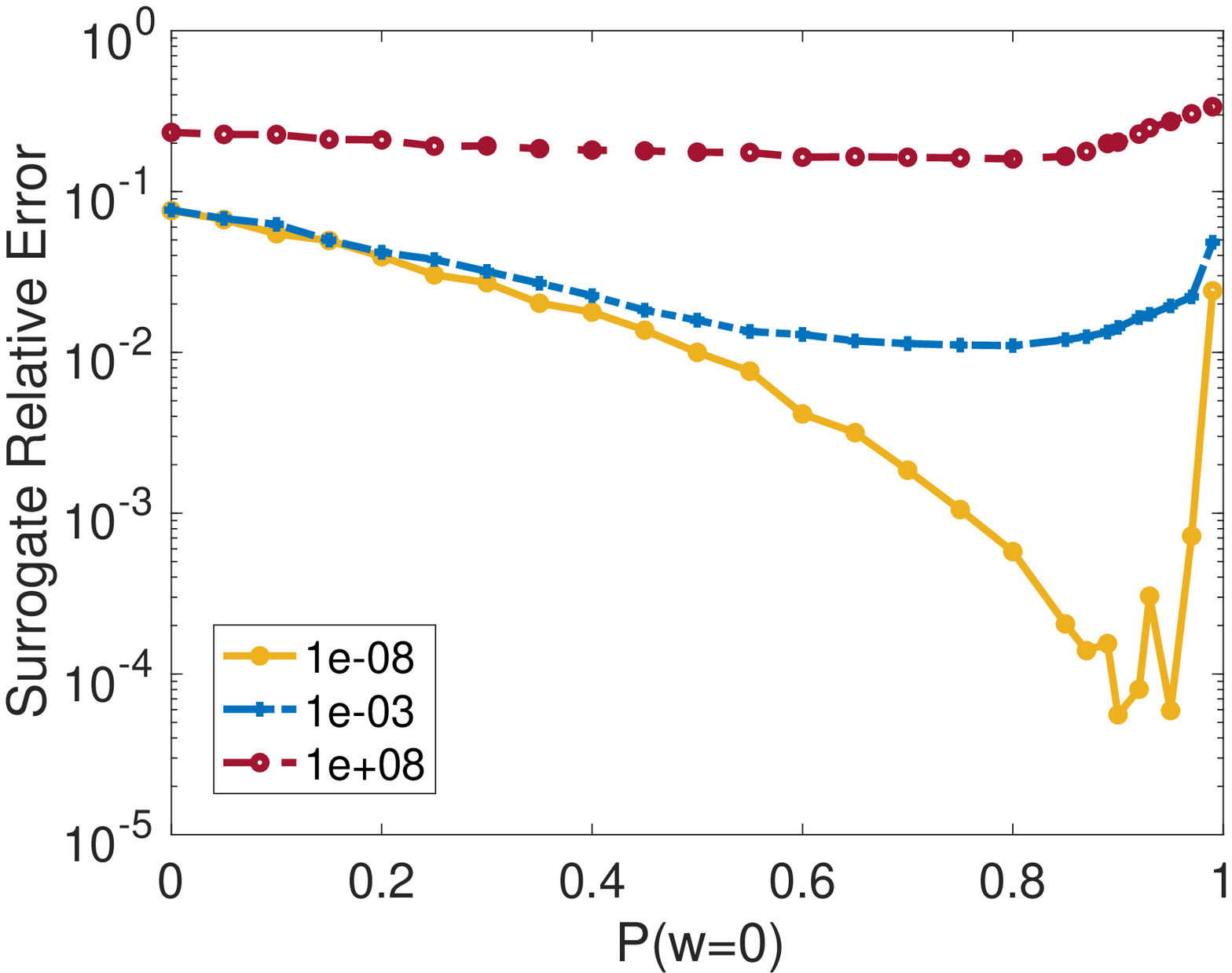}}
   \subfigure{\includegraphics[width=0.49\textwidth]{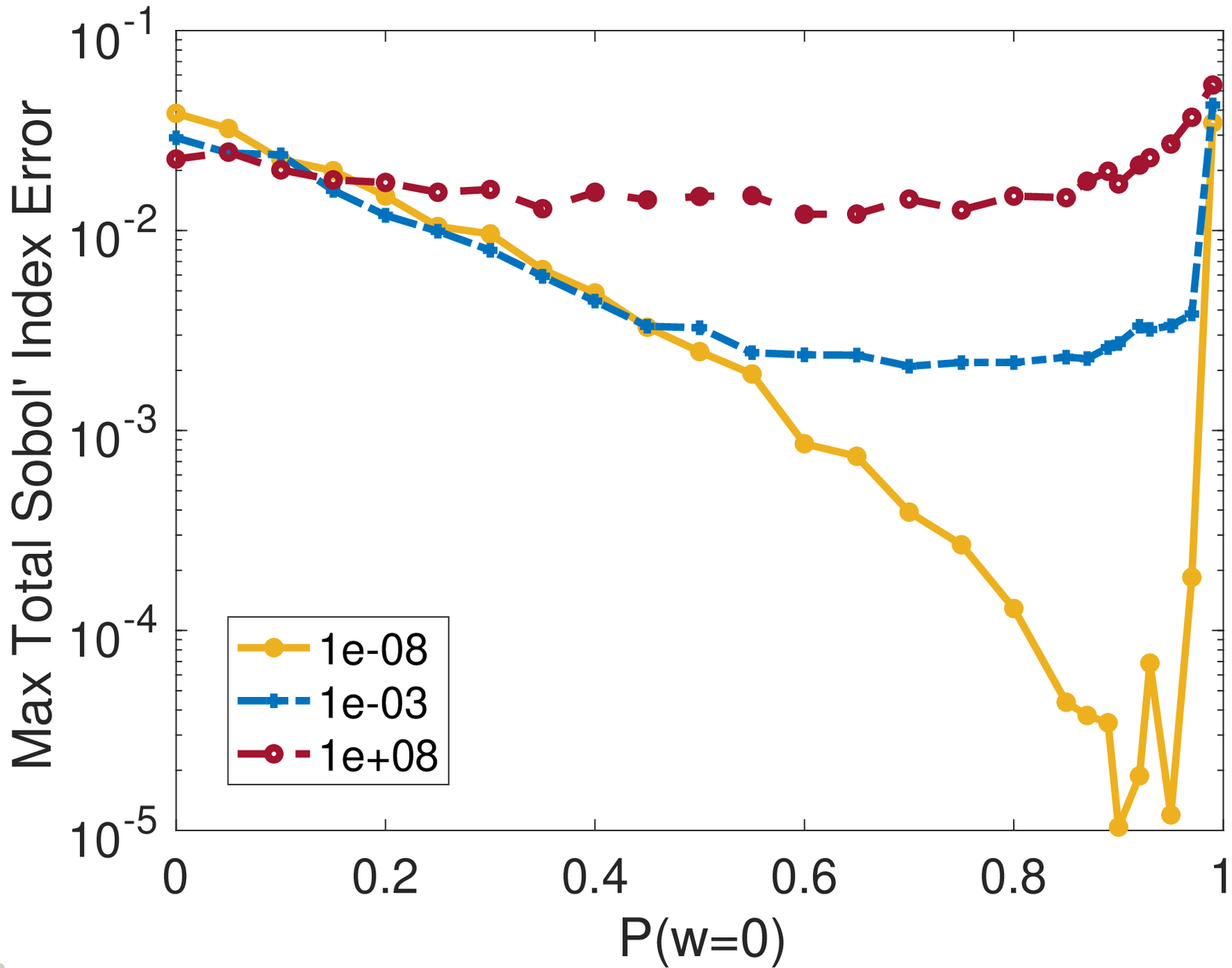}}
 \caption{Example~\eqref{equ:bnch_int} with $d=15$, with low variable interactions  $\delta=10^{-8}$, limited variable interactions $10^{-3}$, and strong variable interactions $10^8$. Surrogate relative error (left) and largest absolute error in total Sobol' index approximation (right).} 
\label{fig:bench_intp}
\end{figure}

We propose sparsifying the hidden layer weight matrix $\mathbf W$ before training the output weights as a means of controlling the influence of variable interactions. Specifically, we replace $\mathbf W$ by a sparse weight matrix $\mathbf{\widetilde W}$ defined as
\begin{equation}\label{equ:sparseW}
\mathbf{\widetilde W} = \mathbf{B}\circ\mathbf{W},
\end{equation}
where $\mathbf{B}$ is a $n\times d$-matrix with independent Bernoulli entries, i.e., 
\begin{equation}\label{equ:sparseB}
B_{ij} = \left\{ \begin{array}{ll} 0 & \mbox{with probability $p$}, \\ 1 & \mbox{with probability $1-p$}\end{array}\right.,
\end{equation}
and where $\circ$ stands for elementwise matrix multiplication. The
sparsification parameter $p\in[0,1)$ determines how sparse the hidden layer
weight matrix is. If $p=0$, then $\mathbf{\widetilde W} = \mathbf{W}$ and
the method reverts to standard ELM; when $p$ is selected to be near $1$, the
weight matrix $\mathbf{\widetilde W}$ is sparse. This technique
strays from the ELM theory which assumes sampling weights from a continuous
probability distribution to guarantee universal approximation.

We implement the sparse weight ELM (SW-ELM) approach, i.e., ELM with
$\mathbf{\widetilde W}$ instead of $\mathbf{W}$ as a hidden layer weight
matrix, on example~\eqref{equ:bnch_int}. \cref{fig:bench_intp} illustrates the
results for a range of values of the sparsification   parameter $p$ and for
three values of $\delta$ corresponding to low variable interactions ($\delta =
10^{-8}$), limited variable interactions ($\delta =10^{-3}$), and strong
variable interactions ($\delta = 10^8$).  In the case of low variable
interactions ($\delta = 10^{-8}$), \cref{fig:bench_intp} shows that using a
sparse weight matrix can dramatically increase the accuracy of both the ELM
surrogate and of the Sobol' indices computed using that surrogate.  Sparsifying
the weight matrix may thus lead  to significant improvements  in GSA results when
dealing with  models with few interaction terms. On the other hand, in cases where interaction
terms are prominent in the model, sparsifying the hidden layer weight matrix
may offer no improvement or may even result in loss of accuracy. 

\subsection{Selection of sparsification parameter}\label{sec:spars}

\cref{fig:bench_intp} suggests a framework for how to implement SW-ELM for a
given model; indeed, for  fixed values of  $\delta$ and varying values of  $p$,
we observe similar trends in the surrogate error and in the total Sobol' index
approximation error. This indicates that we may use the surrogate error as a
guideline when searching for the best sparsification parameter. 

We propose creating a validation set alongside the training set when using
SW-ELM surrogates for GSA.  Subsequently, we select $r$ values of the
sparsification parameter $p$ and construct corresponding sparse weight
matrices.  Then, ELMs are trained using each weight matrix and the validation
set is used to compute the surrogate error for each ELM. This serves as a
diagnostic tool to aid in deciding whether sparsification may improve Sobol'
index approximation. If sparsification yields an improved surrogate error, then
we use the SW-ELM with the lowest surrogate error to approximate Sobol'
indices. If, on the other hand, sparsification does not yield a notable
improvement then we may default to using standard ELM. We summarize the method
in~\cref{alg:spelm}.

\begin{algorithm}[h]
\caption{GSA with sparse weight ELM}\label{alg:spelm}
\textbf{Input:} (i) Model $f$; 
(ii) training set $\{\vec{x}_i,y_i\}_{i=1}^m$; 
(iii) validation set $\{\vec{x}_j',y_j'\}_{j=1}^s$; 
(iv) number of neurons $n$; 
(v) candidate sparsification values $\{p_l\}_{l=1}^r$ (with $p_1=0$)

\textbf{Output:} 
(i) First order Sobol' indices $\{S_k\}_{k=1}^n$ and (ii) total Sobol' indices $\{S_k^\mathrm{tot}\}_{k=1}^n$
\begin{algorithmic}[1]
\State Generate weight matrix $\mat{W}_0$ and bias vector $\vec b$ using a standard normal distribution
\For{$l=1,\dots,r$}
	\State Construct $\mat{W}_l = \mat{B}\circ\mat{W}_0$ with $\mat{B}$ from~\eqref{equ:sparseB} with  $p= p_l$
	\State Determine regularization parameter $\alpha_l$ by the L-curve method
	\State Find output weights $\vec\beta_l$ by training ELM using $\mat{W}_l,\vec{b},\alpha_l$
        (see~\eqref{equ:rlsq})
	\State Compute relative surrogate error $E_l=E_\mathrm{surr}$ (see~\eqref{equ:Esurr}) 
\EndFor
\State Select hidden layer weight matrix $\mat{W}$ and output weights $\vec\beta$ corresponding to sparsification parameter that gives smallest relative surrogate error
\State With $\mat{W},\vec\beta$, and $\vec b$, compute first order Sobol' indices $\{S_k\}_{k=1}^n$ and 
total Sobol' indices $\{S_k^\mathrm{tot}\}_{k=1}^n$ using~\eqref{equ:reg} and~\eqref{equ:tot}, respectively
\end{algorithmic}
\end{algorithm}

\begin{figure}[!]
\centering
\includegraphics[width=1\textwidth]{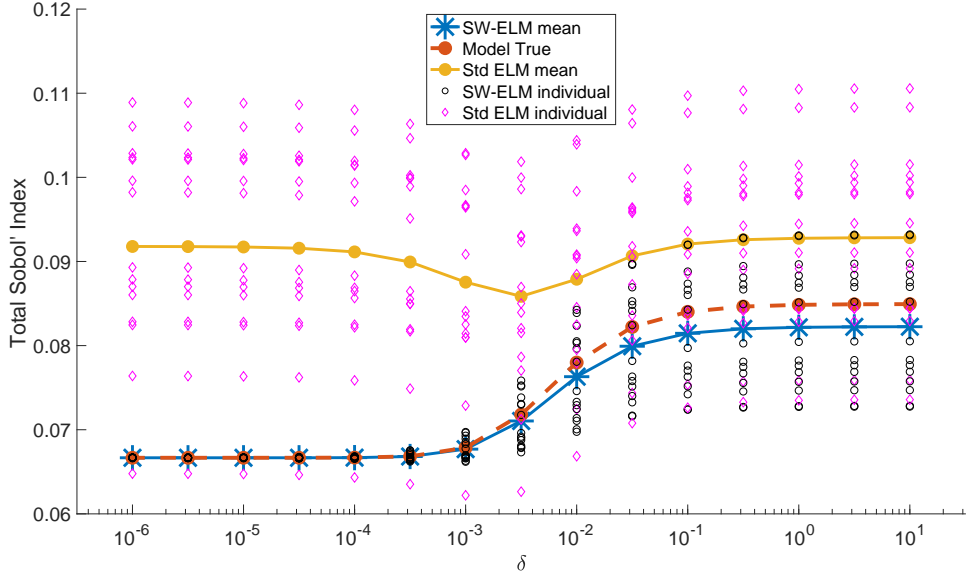}
\caption{Example~\eqref{equ:bnch_int} true analytic total Sobol' indices compared to those approximated by standard ELM and sparse weight ELM (SW-ELM). All approximated indices are plotted along with means for each $\delta$.}
\label{fig:bench_int_fo}
\end{figure}

We illustrate the power of this method by repeating the experiment in~\cref{fig:bench_int} for~\eqref{equ:bnch_int} using SW-ELM instead of ELM. \cref{fig:bench_int_fo} displays the total Sobol' indices for both surrogates as well as the exact indices. 
SW-ELM displays increased accuracy and reduced variance across all tested ranges of interaction strength. 

In the context of GSA, SW-ELM displays better  accuracy and flexibility  than  standard ELM. The proposed implementation is only marginally more costly than standard ELM as there is no need to re-sample training or validation points for each ELM.

\section{Computational results}\label{sec:numerics}

\subsection{Analytic example: Sobol' g-function}

The Sobol' g-function~\cite{SalSob95} 
\begin{equation}\label{equ:gfcn}
f(\vec{x})=\prod_{i=1}^d g_i(x_i) \quad\vec{x}\in[0,1]^d,\quad   g_i(x_i) = \frac{|4x_i-2|+a_i}{1+a_i}, \quad i=1,
\ldots, d,
\end{equation}
is commonly used as a benchmark  to test new methods;   its nonlinearity and
lack of smoothness make approximating its Sobol' indices a significant
challenge. Moreover, we can compare to the true Sobol' indices. Formulas for these can be found in the appendix to~\cite{Saltelli10}. The constants $a_i$, chosen from the interval $(-1,\infty)$, can be
tuned to determine which input variables are important. The closer  $a_i$ is to
$-1$, the more ``important"  $x_i$ becomes. On the other hand, the relative
importance of $x_i$ diminishes with larger values of $a_i$. For our test, we
take the input dimension $d=8$ and let $\vec{a}=$[1, 2, 5, 10, 20, 50, 100,
500] as in~\cite{Sudret08}.
% follow the recommendation of~\cite{Crestaux09}, letting $a_i=\frac{i-2}{2}$
% for $i=1,\dots,d$.

\begin{figure}[h!!]
    \centering
%\subfigure{\includegraphics[width=0.4\textwidth]{./figures/gfunction_lcurve}}
  \includegraphics[width=0.45\textwidth]{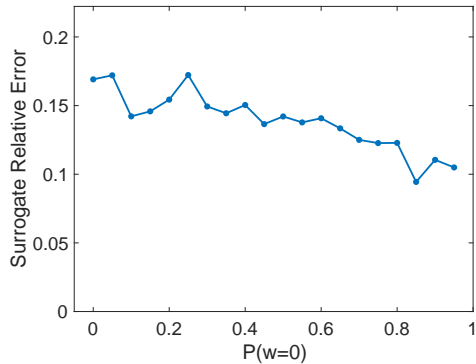}
    \caption{Sparsification test for 8-dimensional g-function~\eqref{equ:gfcn} with $\vec{a}=$[1, 2, 5, 10, 20, 50, 100, 500]. Relative surrogate error is estimated using validation set with 100 points.}
    \label{fig:gfun_tests}
\end{figure}

We train our SW-ELM surrogates with 160 hidden layer neurons using 400 training points and 100 validation points, all sampled by Latin hypercube sampling (LHS). The regularization parameter $\alpha = 10^{-3}$ is selected by the L-curve method. The sparsification parameter $p$ is determined through~\cref{alg:spelm}.  While~\cref{fig:gfun_tests} does not show drastic improvements in error, we nevertheless conclude from this test that sparsifying may improve our approximations when performing GSA; we select the sparsification parameter $p=0.85$. Sobol' indices estimated by SW-ELM are compared to their true analytic values.

\begin{figure}[h!!]
    \centering
\subfigure{\includegraphics[width=0.4\textwidth]{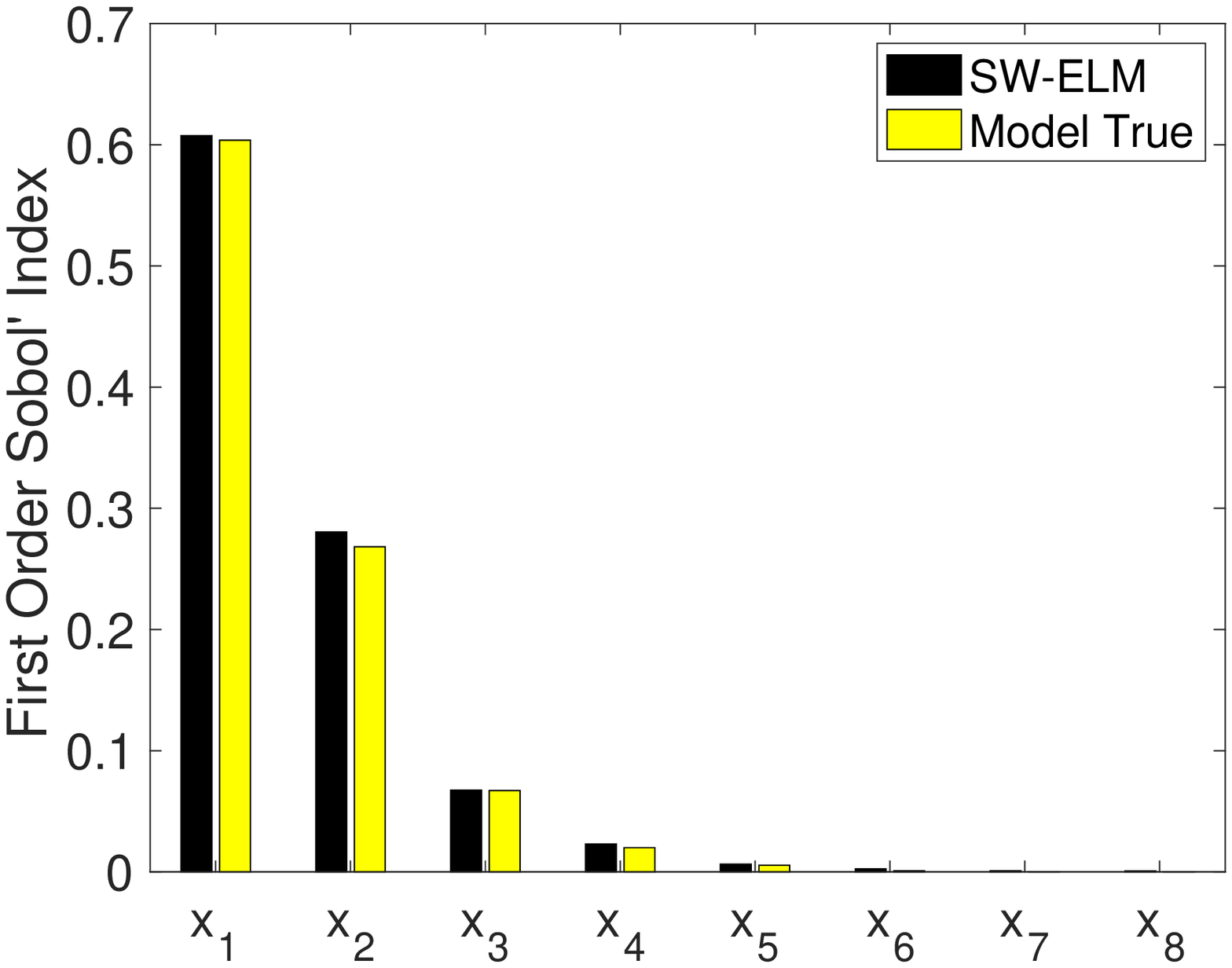}}
    \subfigure{\includegraphics[width=0.4\textwidth]{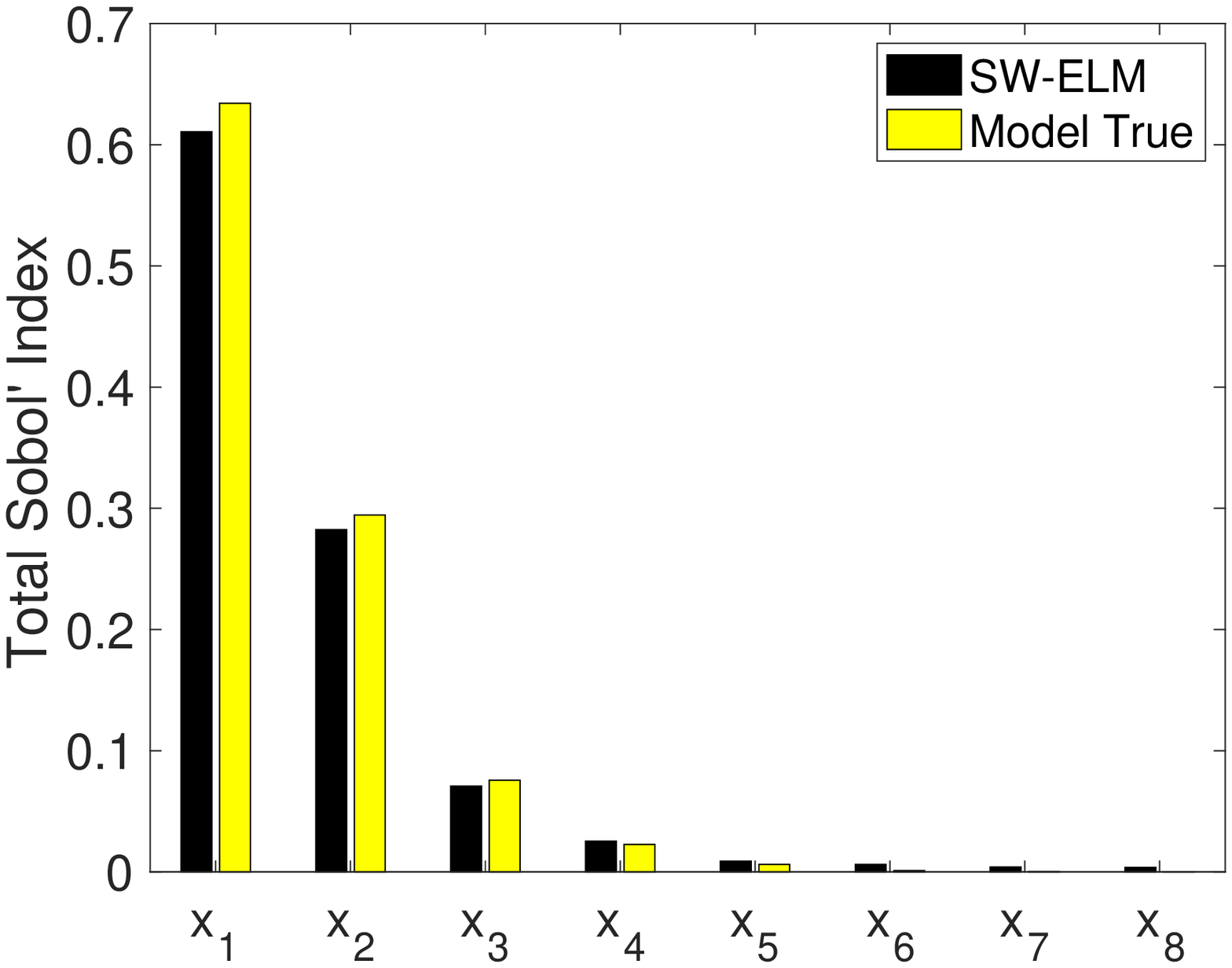}} 
    \caption{8-dimensional g-function~\eqref{equ:gfcn} with $\vec{a}=$[1, 2, 5, 10, 20, 50, 100, 500]. First-order (left) and total (right) Sobol' indices approximated via SW-ELM and computed analytically.}
    \label{fig:gfun_inds}
\end{figure}

As can seen from ~\cref{fig:gfun_inds}, SW-ELM  correctly ranks the Sobol' indices and successfully identifies the most influential input variables.
 For important input variables, each approximated first-order Sobol' index is within $5\%$ relative error of the respective true first-order Sobol' index while each approximated total  Sobol' index is within $7\%$ relative error of the respective true total Sobol' index. 

\begin{figure}[h!!]
    \centering
\subfigure{\includegraphics[width=0.4\textwidth]{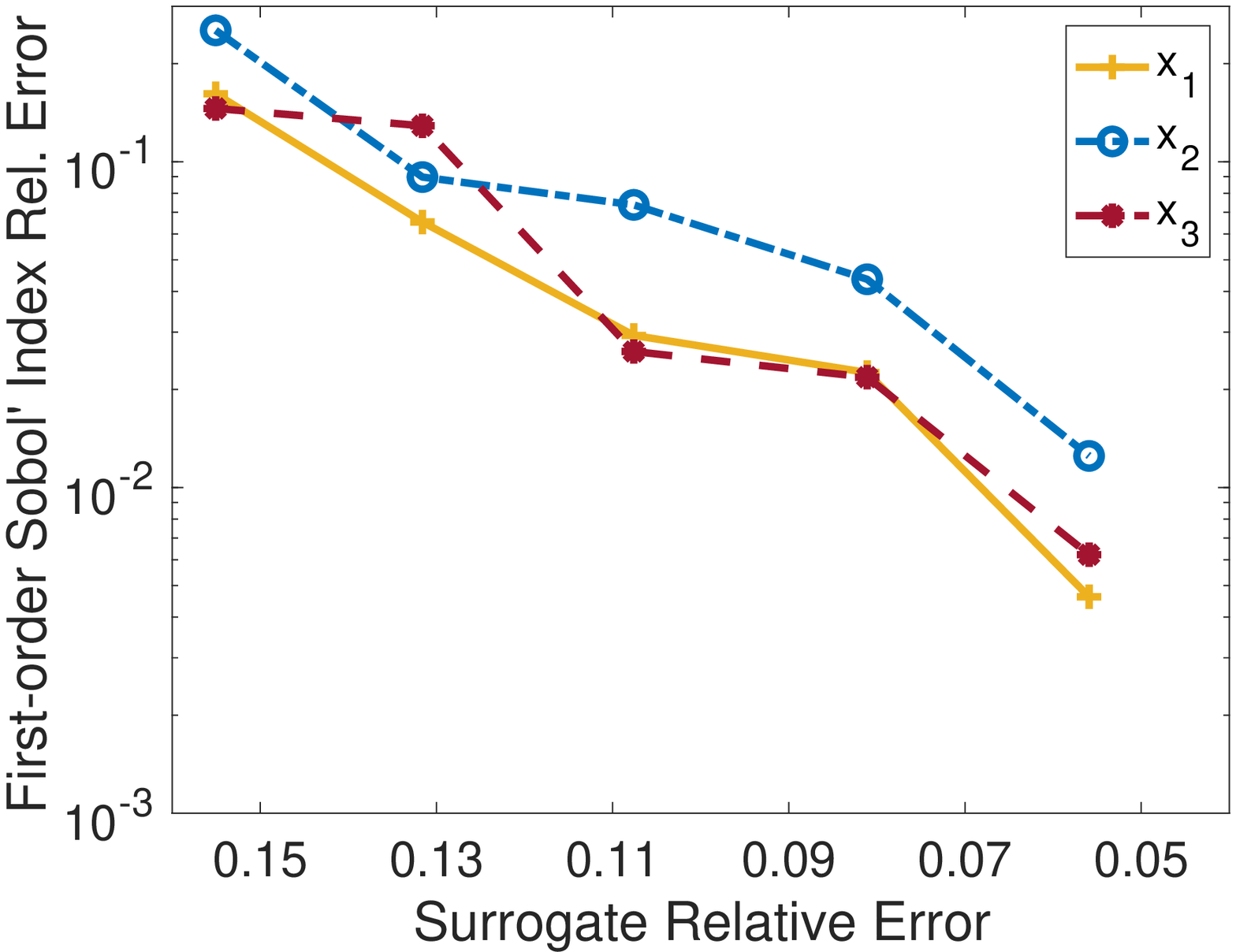}}
    \subfigure{\includegraphics[width=0.4\textwidth]{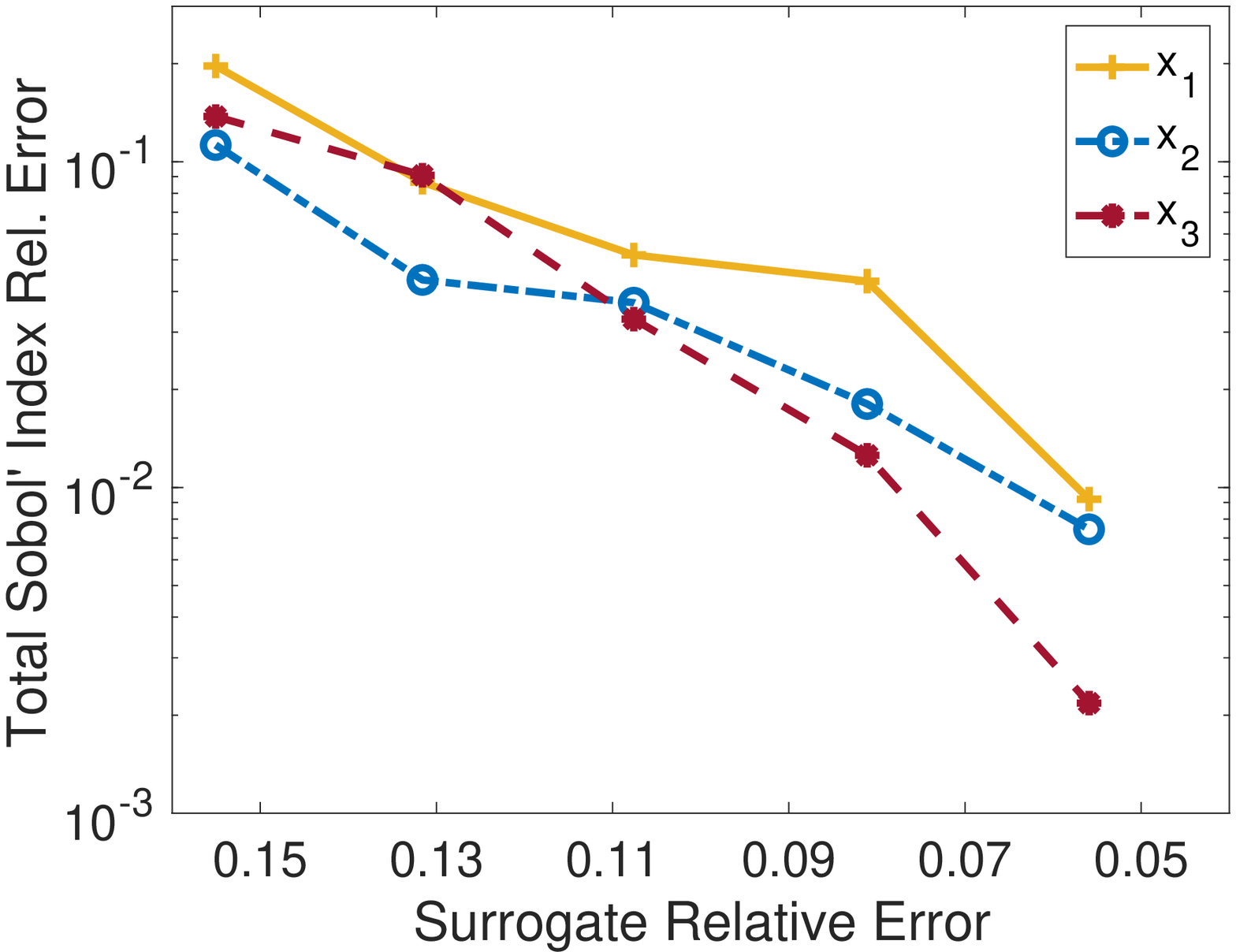}} 
    \caption{8-dimensional g-function~\eqref{equ:gfcn} with $\vec{a}=$[1, 2, 5, 10, 20, 50, 100, 500]. Effect of reducing numerical error on the relative error of first-order (left) and total (right) Sobol' indices estimated by SW-ELM.}
    \label{fig:gfun_ind_err}
\end{figure}

We use this benchmark example to study how numerical accuracy of SW-ELM affects the accuracy of GSA. In the experiment, we train SW-ELMs to a certain accuracy level by adding neurons and training points until the surrogate reaches the accuracy level. Each surrogate uses two training points for every neuron. We decrease the surrogate relative error from $15\%$ to $5\%$ and record the relative error of the first-order and  total Sobol' index estimates for the inputs $x_1,x_2,x_3$. The relations between numerical error and GSA error are shown in~\cref{fig:gfun_ind_err}.
The experiment shows that SW-ELM can provide accurate GSA estimates, especially for the more influential inputs, despite $11\%$ numerical error. Still, reducing the numerical error has a demonstrative effect on GSA. As we decrease the error to $5\%$, the Sobol' indices become much more accurate.

%%%%%%%%%%%%%%%%%%%%%%%%%%
%%%%%%%%%%%%%%%%%%%%%%%%%%

\subsection{Genetic oscillator}
We examine the performance of the proposed approach on a
challenging application problem from biochemistry, specifically, 
the genetic oscillator system which  describes the time
evolution of molecular species involved in the regulation of circadian
rhythm~\cite{Vilar02}. The reaction network consists of sixteen reactions and
involves nine species. The reactions, the corresponding propensity functions,
and the nominal values of the rate parameters are given in~\cref{tab:genetic}.
\begin{table*}[h!]\label{tab:genetic}
\centering
 \begin{tabular}[!htb]{cc|cc}
 \hline
Reaction & Propensity Function & Parameter & Value \\ [0.5ex] 
 \hline
$P_a \to P_a + mRNA_a$ & $\alpha_AP_a$  & $\alpha_A$ & $50.0$\\ 
$P_r \to P_r + mRNA_r$ & $\alpha_R P_r$   & $\alpha_R$ & $0.01$\\
 $mRNA_a \to mRNA_a + A$  & $\beta_A mRNA_a$   & $\beta_A$ & $50.0$ \\
$mRNA_r \to mRNA_r + R$  & $\beta_R mRNA_r$  & $\beta_R$ & $5.0$\\
 $A + R \to C$  & $\gamma_C AR$ & $\gamma_C$ & $20.0$ \\
$P_a + A \to P_{a-}A$  & $\gamma_A P_a A$ & $\gamma_A$ & $1.0$ \\
 $P_{a-}A \to P_a + A$  & $\theta_A P_{a-}A$ & $\theta_A$ & $50.0$  \\
 $P_r + A \to P_{r-}A$  & $\gamma_R P_r A$  & $\gamma_R$ & $1.0$ \\
$P_{r-}A \to P_r + A$  & $\theta_R P_{r-}A$  & $\theta_R$ & $1.0$\\
 $A \to \emptyset$  & $\delta_A A$  & $\delta_A$ & $1.0$ \\
 $R \to \emptyset$  & $\delta_R R$  & $\delta_R$ & $0.2$  \\
$mRNA_a \to \emptyset$  & $\delta_{MA} mRNA_a$ &  $\delta_{MA}$ & $10.0$ \\
$mRNA_r \to \emptyset$  & $\delta_{MR} mRNA_r$  & $\delta_{MR}$ & $0.5$ \\
  $C \to R$  & $\delta_A' C$ & $\delta_A'$ & $1.0$ \\
 $P_{a-}A \to P_{a-}A + mRNA_a$ & $\alpha_a\alpha_A P_{a-}A$  &  $\alpha_a$ & $10.0$\\
 $P_{r-}A \to P_{r-}A + mRNA_r$  & $\alpha_r\alpha_R P_{r-}A$  & $\alpha_r$ & $5000$ \\ [1ex] 
 \hline
\end{tabular} 

\vspace{1mm}
\caption{Genetic oscillator reactions, propensity functions, parameters and nominal values of the  parameters~\cite{Sheppard12,Merritt21}.}
\end{table*}
We consider the reaction rate equations (RREs), described by a nonlinear system of
ordinary differential equations (ODEs), for this reaction network; see
e.g.,~\cite{Vilar02},~\cite{Sheppard12}, or~\cite{Merritt21}, which we follow in the specific
problem formulation used in the present study.  We focus on the uncertainty in
the reaction rate parameters.  A uniform distribution is attached to each rate
parameter on an interval given by a $\pm 5\%$ perturbation from the corresponding nominal value.  
We use a random vector
$\vec{x} \in \R^{16}$, whose entries are uniformly distributed on $[0, 1]$, to
parameterize the uncertainty in the reaction rates.  The vector of reaction
rates is obtained by applying a linear transformation to $\vec{x}$ that maps
the entries of $\vec{x}$ to the respective physical ranges.
In the present study, we consider the quantity of interest (QoI) given by 
\begin{equation}\label{equ:qoi}
f(\vec{x})=\frac{1}{T}\int_0^TR(t;\vec{x})\,dt,
\end{equation}
where $R(t;\vec{x})$ is the concentration of the species $R$ present in
the system at time $t$, and $T$ is the final simulation time. 
Notice that computing $R(t; \vec x)$, for $t \in [0, T]$, requires solving
the system of RREs with the reaction rates set according to $\vec x$. 
This system, for the present application, is
given by a stiff system of nonlinear ODEs.
Thus, evaluations of the QoI in~\eqref{equ:qoi} are   
computationally expensive.

\begin{figure}[h!!]
    \centering
\includegraphics[width=0.49\textwidth]{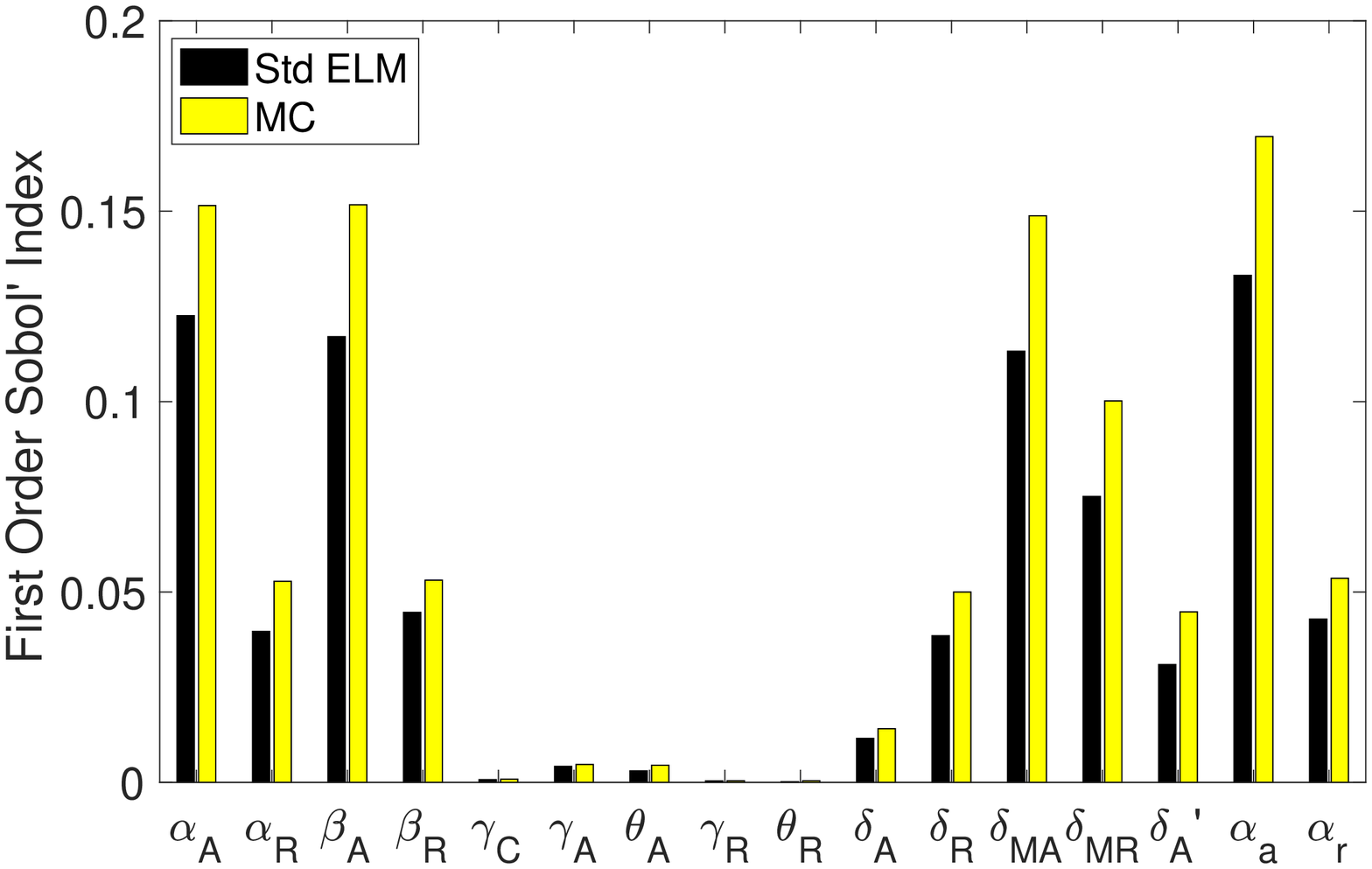}
\includegraphics[width=0.49\textwidth]{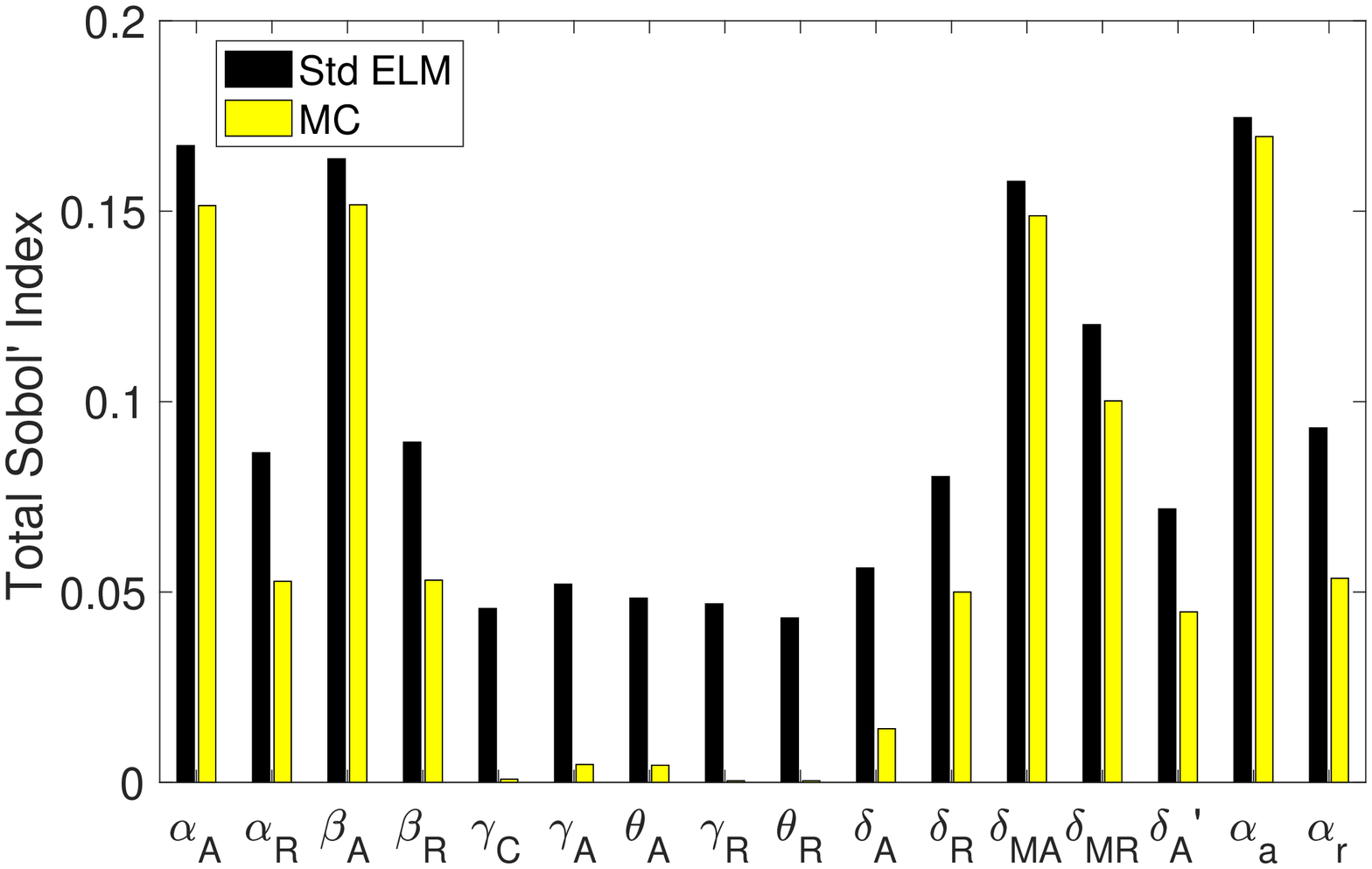} \\
\vskip 1mm
\hskip 0.65mm\includegraphics[width=0.49\textwidth]{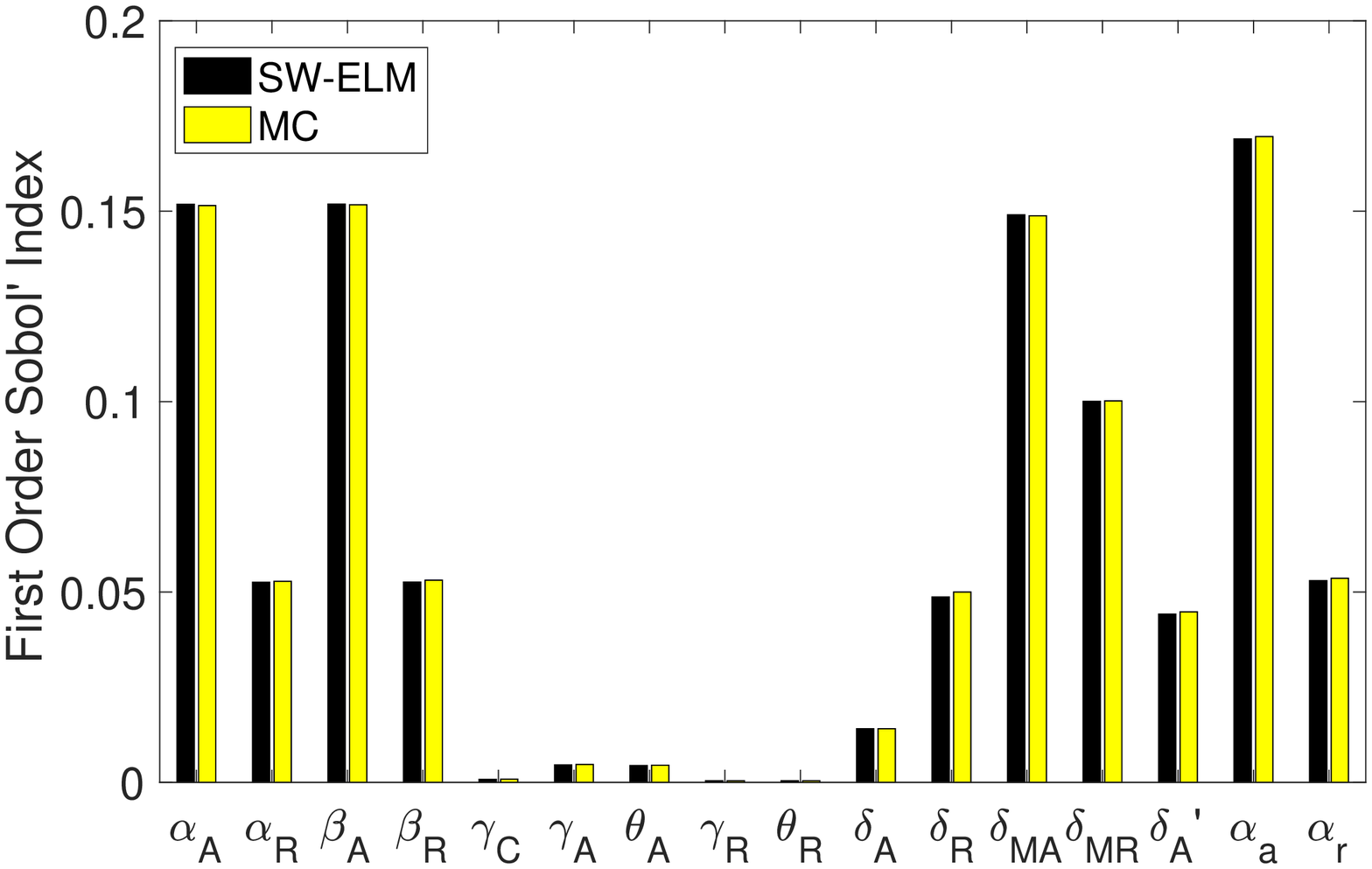}
\includegraphics[width=0.49\textwidth]{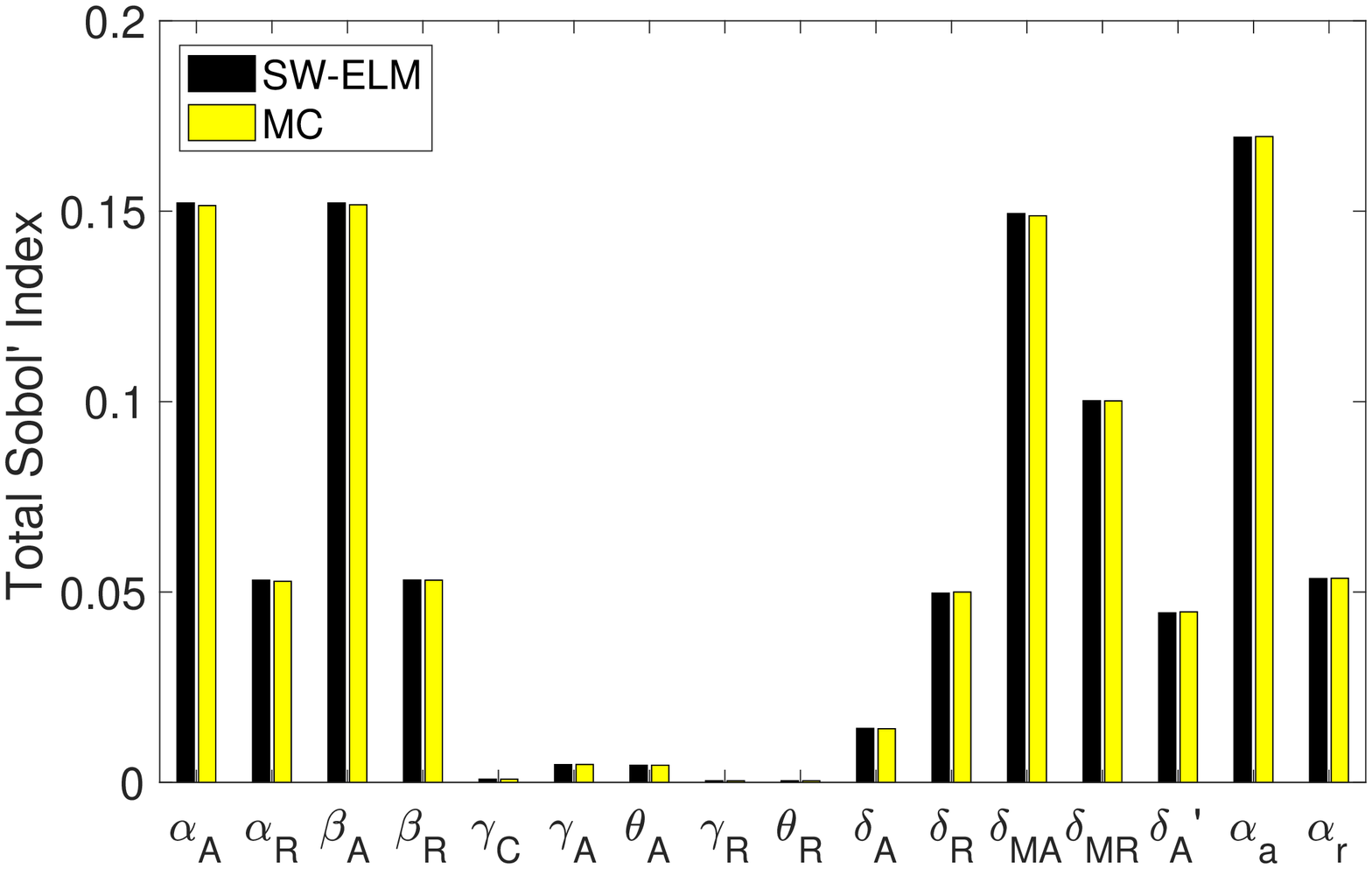} \\
\vskip 1mm
\hskip 0.65mm\includegraphics[width=0.49\textwidth]{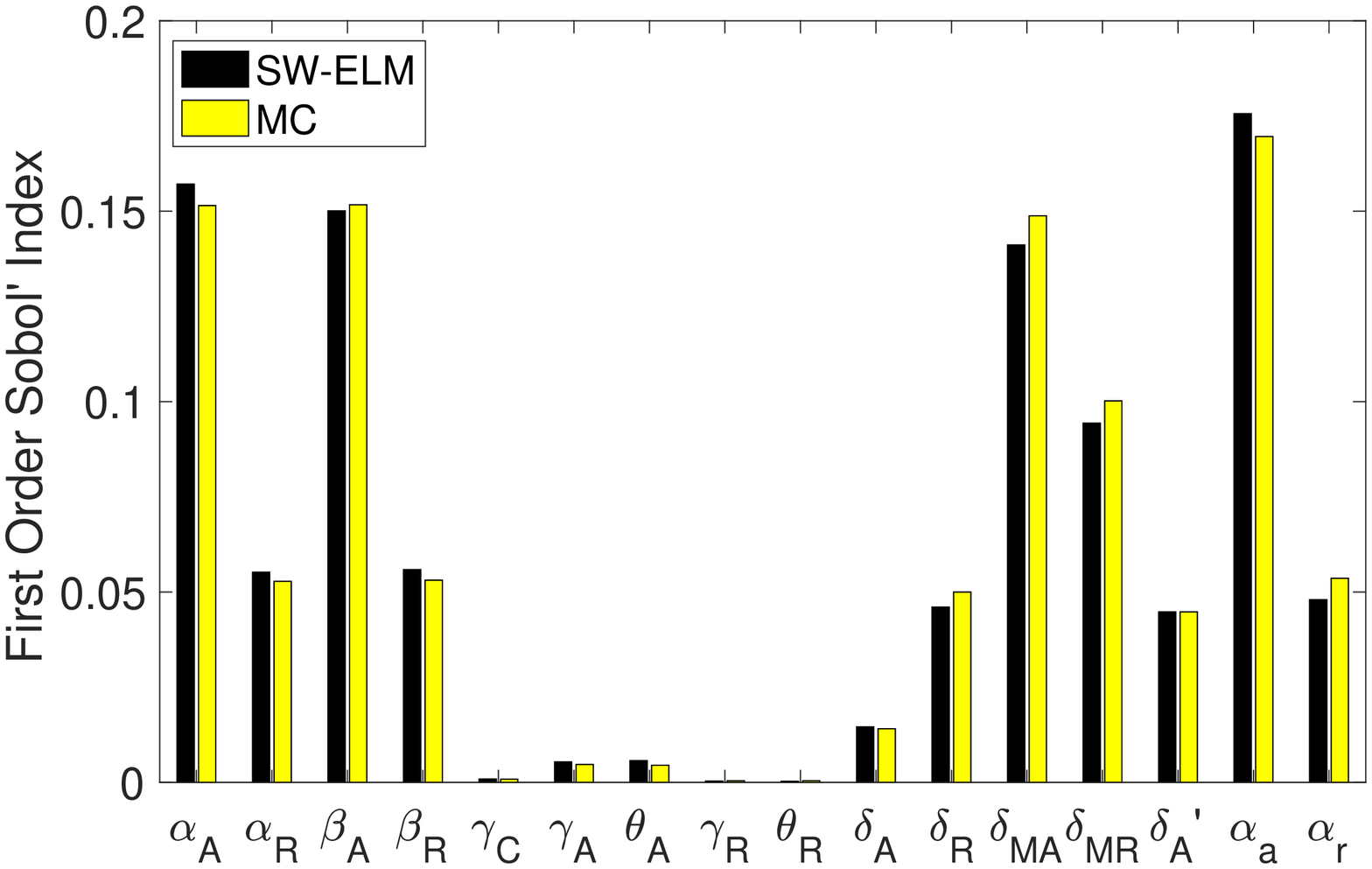}
\includegraphics[width=0.49\textwidth]{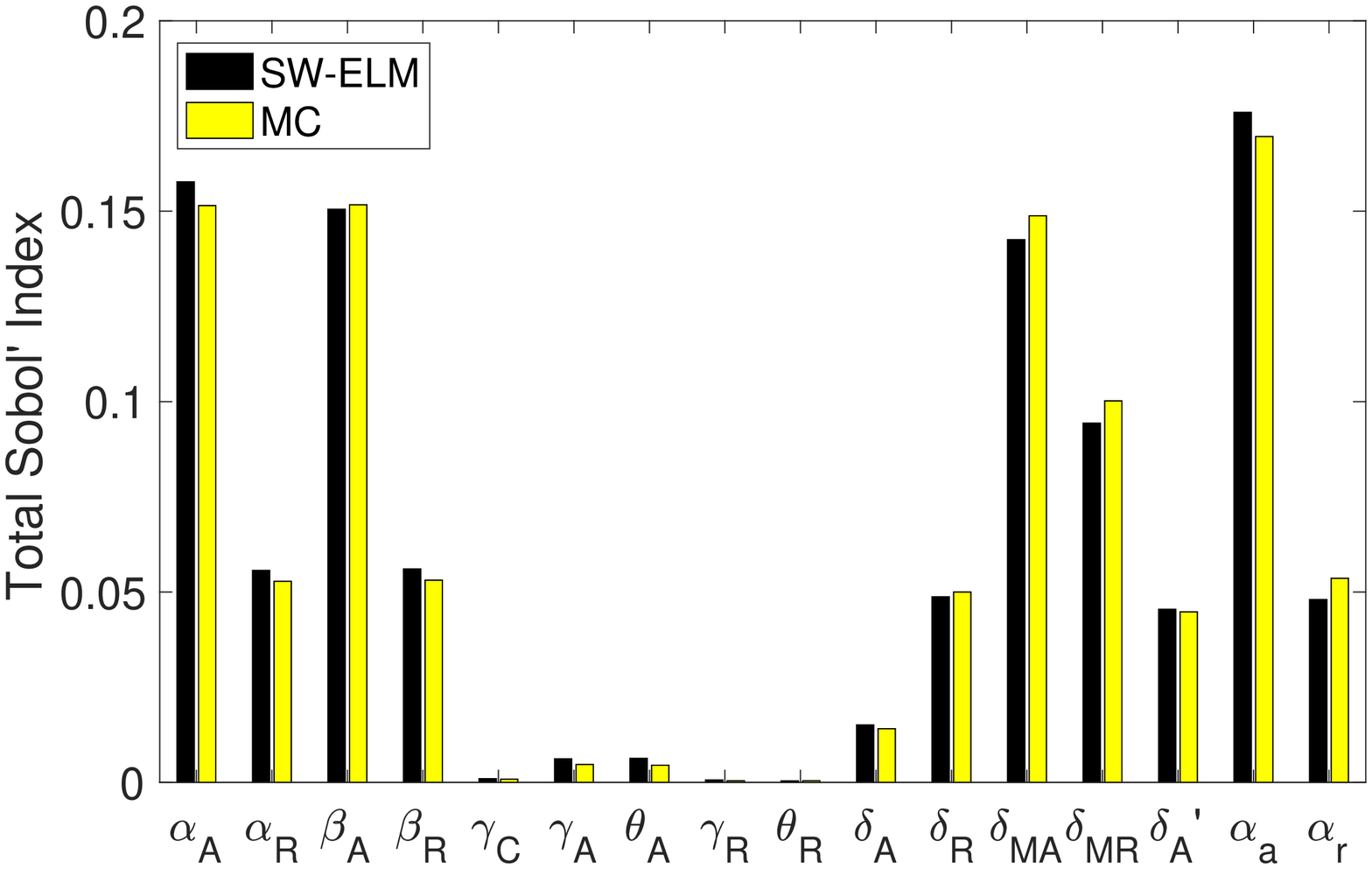}
 \caption{Top row: first-order (left) and total Sobol's (right) indices corresponding to the QoI~\eqref{equ:qoi} approximated via standard ELM, using training set with 3000 points, along with those approximated via Monte Carlo using $10^6$ sample points.
   Middle row: same experiment with SW-ELM, with 100 point validation set, instead of standard ELM.
  Bottom row: same experiment with SW-ELM instead of ELM and training size of 150 instead of 3000.}
    \label{fig:go}
\end{figure}

\cref{fig:go}~(top row) presents approximations for first-order and total
Sobol' indices when using standard ELM. The experiment uses 3000 training
points and 1000 neurons. The regularization parameter $\alpha=10^{-4}$ is
selected by the L-curve method. We observe that approximation using standard
ELM leads to overestimates, particularly for unimportant input variables, of
the total Sobol' indices and underestimates of the first-order Sobol'
indices; this is a clear case of standard ELM overestimating contributions to
the output variance by higher order interactions.

\begin{figure}[h!!]
    \centering
    \includegraphics[width=0.49\textwidth]{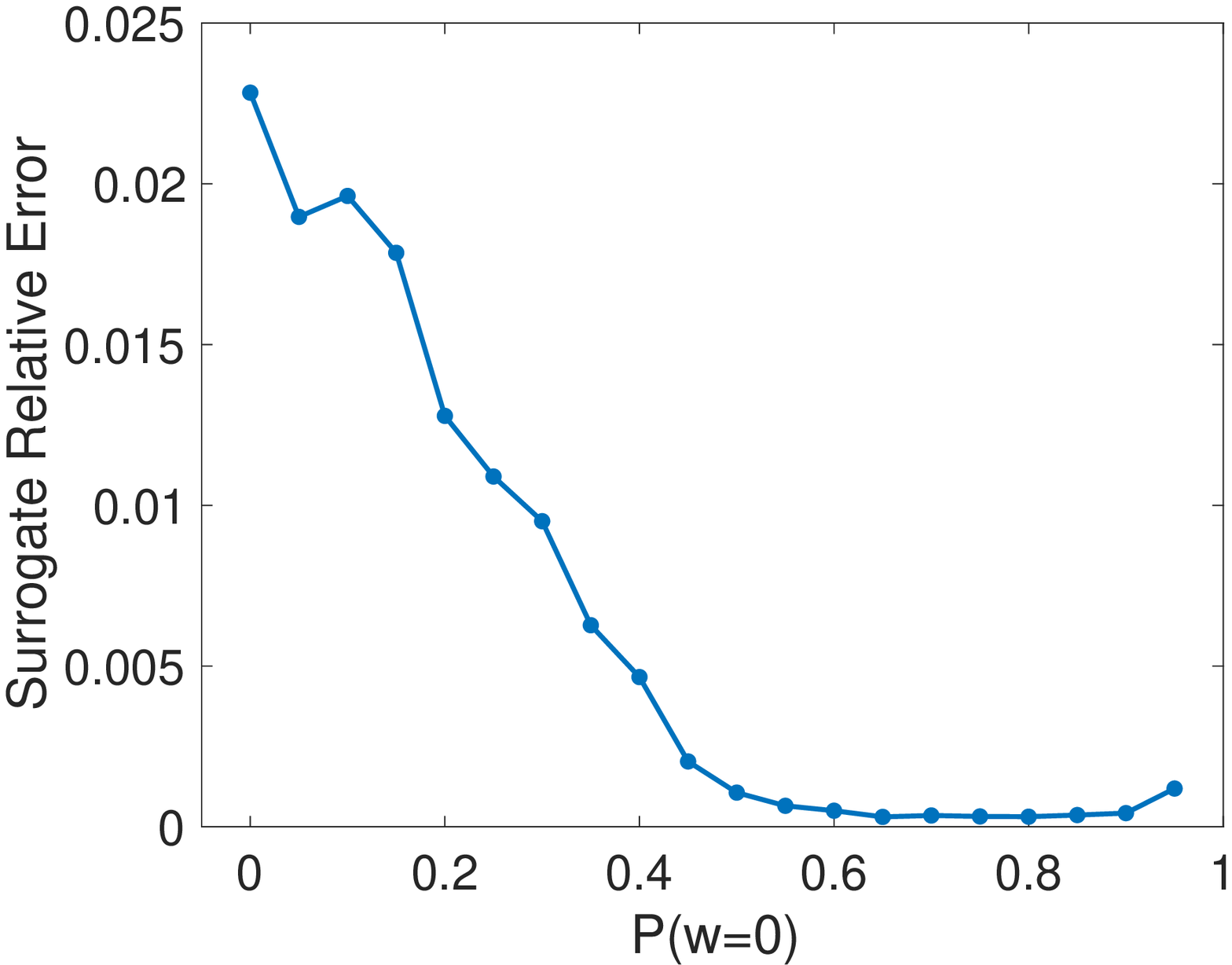}
	\includegraphics[width=0.49\textwidth]{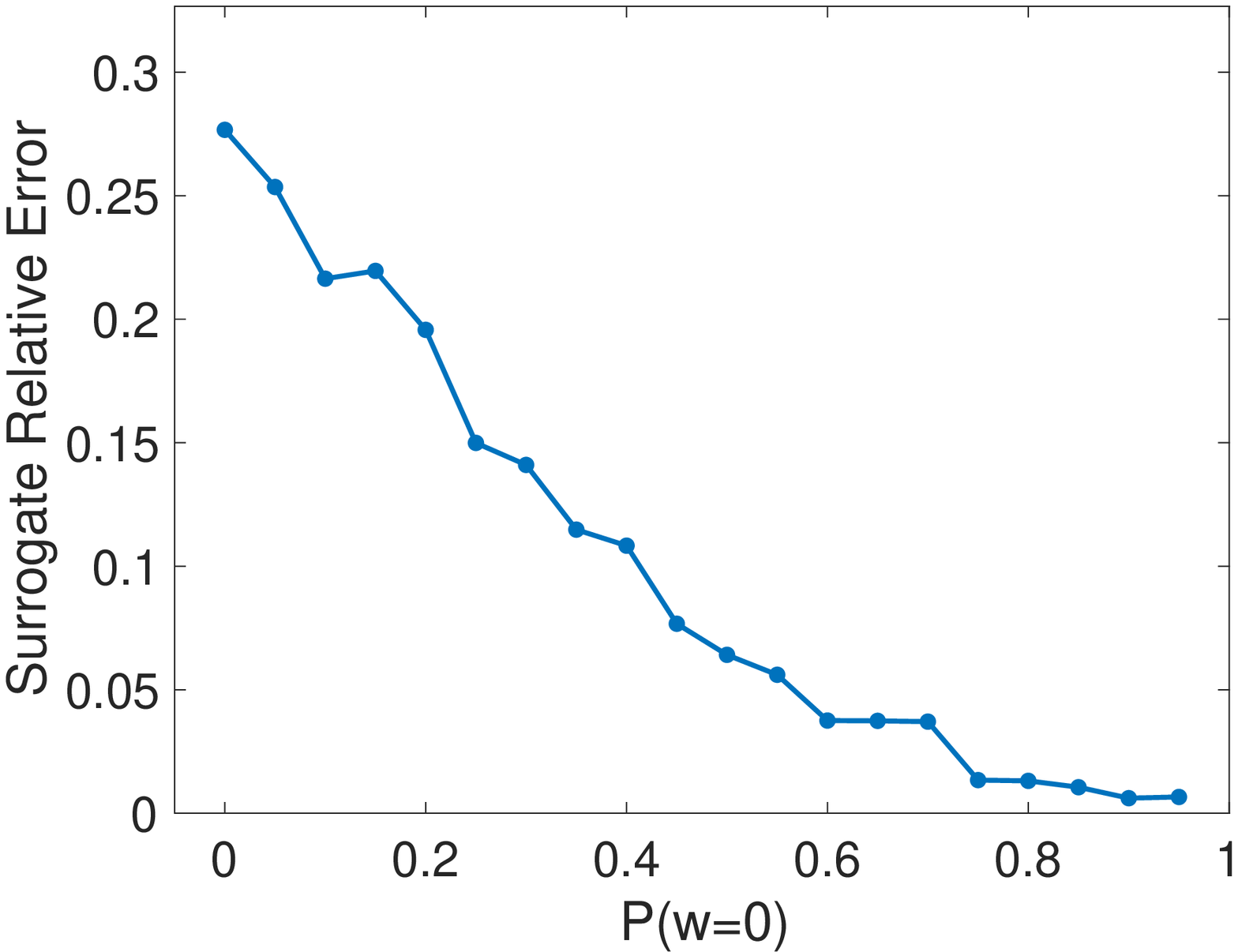}
    \caption{Sparsification test for 16-dimensional QoI~\eqref{equ:qoi}. Left: ELMs trained with training size 3000. Right: ELMs trained with training size 150. Relative surrogate error is is estimating using validation set with 100 points.}
    \label{fig:go_tests}
\end{figure}

Let us instead tackle this problem using SW-ELM; we train our surrogates with
again 3000 training points, using 1000 hidden layer neurons, and 100 validation
points, all sampled by LHS. Regarding sparsification, we
follow~\cref{alg:spelm}. The results of these tests are displayed
in~\cref{fig:go_tests}~(left), and show that as the hidden layer weight
matrix become more sparse, the error in the surrogate improves dramatically; we
select the sparsification parameter $p=0.8$.  \cref{fig:go}~(middle row),
displays the resulting  approximated Sobol' indices  and compares with the
Sobol' indices approximated by Monte Carlo methods with $10^6$ sample points.
Given that the Sobol' indices computed via Monte Carlo methods are a close
representation of the true Sobol' indices, SW-ELM provides an accurate
approximation  of the indices. 
 
We can further challenge the performance of SW-ELM. We train surrogates with
with only 150 training points, using 50 hidden layer neurons, and 100
validation points. The results of sparsification tests are displayed
in~\cref{fig:go_tests}~(right). We select the sparsification parameter $p=0.9$.
\cref{fig:go}~(bottom row) displays the resulting  approximated Sobol' indices.
SW-ELM with 150 training points or 3000 training points gives very similar
results.  This indicates that SW-ELM is extremely efficient for this 
problem---a modest number of training samples are sufficient for obtaining
accurate sensitivity analysis results.

As shown in~\cref{fig:go}, the first-order Sobol' indices and the total Sobol'
indices are very close to each other for this appplication. This suggests that
there are few, if any, higher order variable interactions present in the QoI
function~\eqref{equ:qoi}. The fact that the sparsity of the weight matrix had a
significant influence over the quality of the ELM surrogate, as seen
in~\cref{fig:go_tests}, further supports this.  Here sparsification is not
simply helpful; it is essential and allows us to drastically decrease the
number of sample points needed to obtain accurate approximations. 

\subsection{High-dimensional example}
\label{sec:highdim}
In this section, we use our proposed approach to peform GSA in a model governed by a
system of ordinary differential equations (ODEs) with has $50$ uncertain
inputs.  Specifically, we consider a linear homogeneous ODE system,
\begin{equation}\label{equ:lode}
\dot{\boldsymbol{x}} = -\mat A \boldsymbol{x}, \quad \boldsymbol{x}_0 = \boldsymbol{x}(0),
\end{equation}
where $\boldsymbol{x}(t) = [x_1(t),\quad\dots,\quad x_{50}(t)]^\top$ and $\mat
A\in\R^{50\times 50}$. Note that the solution of this is system is given by
$\boldsymbol{x}(t) = e^{-t\mat A}\boldsymbol{x}_0$.

In the present example, we assume $\mat{A}$ is symmetric. Thus, $\mat A$ has 
spectral decomposition $\mat A  = \mat Q\Lambda\mat Q^\top$, where $\mat Q$ is orthogonal and $\Lambda$ is diagonal with 
the eigenvalues $\lambda_1,\dots,\lambda_{50}$ of $\mat A$ as its diagonal entries. 
This allows us to express the solution as
\begin{equation}\label{equ:eigsol}
\boldsymbol{x}(t) = \mat Q e^{-t\boldsymbol\Lambda}\mat Q^\top\boldsymbol{x}_0,
\end{equation}
where $e^{-t\boldsymbol\Lambda}$ is the diagonal matrix taking $e^{-t\lambda_1},...,e^{-t\lambda_{50}}$ as its diagonal entries. 
In this example, we consider a symmetric matrix $\mat A$ with a known 
orthogonal matrix $\mat Q$ of eigenvectors and uncertain eigenvalues. 
That is, the diagonal entries of $\mat\Lambda$ in~\cref{equ:eigsol}, are uncertain.
We assume these eigenvalues take nominal values $\lambda_k =\frac{1}{k},\quad k=1,\dots,50$. 
Each eigenvalue is assumed to follow a uniform distribution on an interval given by a $\pm 5\%$ perturbation from 
the respective nominal value. We consider the quantity of interest,
\begin{equation}\label{equ:highdimqoi}
f(\boldsymbol{\theta}) = x_{50}(10),
\end{equation}
which is the last entry in the solution vector 
$\boldsymbol{x}(t)=\mat Q e^{-t\boldsymbol\Lambda}\mat Q^\top\boldsymbol{x}_0$
at time $t=10$. For further details on the setup of the problem,
see~\cref{sec:highdim_details}.

\begin{figure}[h!!]
    \centering
    \includegraphics[width=0.5\textwidth]{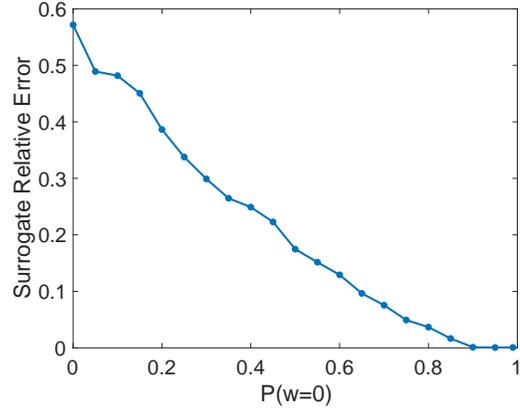}
        \caption{Sparsification test for~\eqref{equ:highdimqoi} on a 100 point validation set. ELMs trained with 700 training points and 350 neurons. Surrogate relative error reaches a minimum of $2.73\times 10^{-4}$ on the validation set at $p=0.95$.}
        \label{fig:eigs_pplot}
\end{figure}

\begin{figure}[h!!]
    \centering 
	\includegraphics[width=0.7\textwidth]{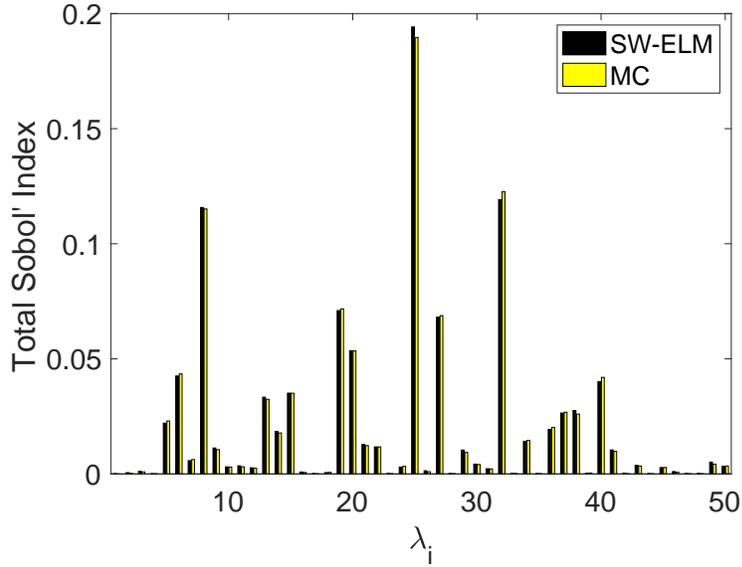}
    \caption{Total Sobol' indices computed for SW-ELM compared to MC estimates. SW-ELM with $p=0.95$ trained on 700 training points with 350 neurons.}
    \label{fig:eigs_tests}
\end{figure}

In~\cref{fig:eigs_tests}, we see the results of the sparsification test. The experiment uses 700 training points, 350 neurons, and 100 validation points to compute surrogate error. As show in~\cref{fig:eigs_pplot}, the surrogate error reaches a minimum when $p=0.95$, indicating that a very sparse weight matrix should be used.
~\cref{fig:eigs_tests} demonstrates the accuracy of the total Sobol' indices estimated using SW-ELM. These are compared to total Sobol' indices estimated using Monte Carlo integration which we treat as the ground truth. The Sobol' indices reveal the absence of variable interactions, justifying the use of a very sparse weight matrix.

\section{Conclusion}
The use of SW-ELM as a surrogate in the context of variance based GSA allows
the user to completely eschew Monte Carlo integration which is a perennial
bottleneck in this field. The numerical results we present above make a strong
case for further study of the approach. In particular, there is in general very
little known theoretically  regarding the link from \eqref{equ:approxf}
(accurate surrogate) to \eqref{equ:approxs} (accurate GSA). In our particular
case, we note further that SW-ELM is not strictly covered by the existing
theory ~\cite{HuangELM} as the distribution we use for $\mat{W}$ is not
continuous. 

A key point in our approach is the ability to analytically compute the Sobol'
indices of $\hat f$. This is where the assumption of independent uniformly
distributed parameters is crucial. Further work will reveal to what extent this
assumption might be weakened. A deeper understanding of the proposed
sparsification process might also result in further improvement of the proposed
approach.




\bibliographystyle{elsarticle-num} 
\bibliography{refs}

\begin{thebibliography}{10}
\expandafter\ifx\csname url\endcsname\relax
  \def\url#1{\texttt{#1}}\fi
\expandafter\ifx\csname urlprefix\endcsname\relax\def\urlprefix{URL }\fi
\expandafter\ifx\csname href\endcsname\relax
  \def\href#1#2{#2} \def\path#1{#1}\fi

\bibitem{Saltelli08}
A.~Saltelli, M.~Ratto, T.~Andres, F.~Campolongo, J.~Cariboni, D.~Gatelli,
  M.~Saisana, S.~Tarantola, Global sensitivity analysis: the primer, John Wiley
  \& Sons, 2008.

\bibitem{IoossLeMaitre15}
B.~Iooss, P.~Le~{M}a{\^i}tre, {A Review on Global Sensitivity Analysis
  Methods}, Springer US, Boston, MA, 2015, pp. 101--122.

\bibitem{IoossSaltelli17}
B.~Iooss, A.~Saltelli, Introduction to sensitivity analysis, in: R.~Ghanem,
  D.~Higdon, H.~Owhadi (Eds.), Handbook of uncertainty quantification,
  Springer, 2017, pp. 1103--1122.

\bibitem{SalSob95}
A.~Saltelli, I.~Sobol', {Sensitivity analysis for nonlinear mathematical
  models: numerical experience}, Matematicheskoe Modelirovanie 7~(11) (1995)
  16--28.

\bibitem{Sobol01}
I.~Sobol', Global sensitivity indices for nonlinear mathematical models and
  their {M}onte {C}arlo estimates, Mathematics and Computers in Simulation
  55~(1--3) (2001) 271--280, the Second IMACS Seminar on Monte Carlo Methods.

\bibitem{PrieurTarantola17}
C.~Prieur, S.~Tarantola, {Variance-based sensitivity analysis: Theory and
  estimation algorithms}, in: R.~Ghanem, D.~Higdon, H.~Owhadi (Eds.), Handbook
  of Uncertainty Quantification, Springer, 2017, pp. 1217--1239.

\bibitem{HartGremaudDavid19}
J.~Hart, P.~Gremaud, T.~David, Global sensitivity analysis of high-dimensional
  neuroscience models: An example of neurovascular coupling, Bulletin of
  mathematical biology 81~(6) (2019) 1805--1828.

\bibitem{Tabandeh22}
A.~Tabandeh, N.~Sharma, P.~Gardoni, Uncertainty propagation in risk and
  resilience analysis of hierarchical systems, Reliability Engineering \&
  System Safety 219 (2022) 108208.

\bibitem{GratietMarelliSudret17}
L.~Le~Gratiet, S.~Marelli, B.~Sudret, {Metamodel-based sensitivity analysis:
  polynomial chaos expansions and Gaussian processes}, in: R.~Ghanem,
  D.~Higdon, H.~Owhadi (Eds.), Handbook of Uncertainty Quantification,
  Springer, 2017.

\bibitem{Sargsyan2017}
K.~Sargsyan, Surrogate models for uncertainty propagation and sensitivity
  analysis, in: R.~Ghanem, D.~Higdon, H.~Owhadi (Eds.), Handbook of uncertainty
  quantification, Springer, 2017.

\bibitem{Sudret08}
B.~Sudret, {Global sensitivity analysis using polynomial chaos expansions},
  journal = {Reliability Engineering \& System Safety} 93~(7) (2008) 964--979,
  bayesian Networks in Dependability.

\bibitem{Crestaux09}
T.~Crestaux, O.~Le~Ma\^{i}tre, J.-M. Martinez, {Polynomial chaos expansion for
  sensitivity analysis}, Reliability Engineering \& System Safety 94~(7) (2009)
  1161--1172, special Issue on Sensitivity Analysis.

\bibitem{friedman91}
J.~H. Friedman, Multivariate adaptive regression splines, The Annals of
  Statistics 19~(1) (1991) 1--141, with discussion and a rejoinder by the
  author.

\bibitem{HartAlexanderianGremaud17}
J.~{Hart}, {A. Alexanderian}, P.~{Gremaud}, {Efficient computation of Sobol'
  indices for stochastic models}, SIAM Journal on Scientific Computing. 39~(4)
  (2017) A1514--A1530.

\bibitem{Marrel2IoossLaurentEtAl09}
A.~Marrel, B.~Iooss, B.~Laurent, O.~Roustant, {Calculations of Sobol' indices
  for the Gaussian process metamodel}, Reliability Engineering \& System Safety
  94~(3) (2009) 742--751.

\bibitem{OakleyOHagan04}
J.~E. Oakley, A.~O'Hagan, {Probabilistic sensitivity analysis of complex
  models: a Bayesian approach}, Journal of the Royal Statistical Society:
  Series B (Statistical Methodology) 66~(3) (2004) 751--769.

\bibitem{JinChenSudjianto04}
R.~Jin, W.~Chen, A.~Sudjianto, Analytical metamodel-based global sensitivity
  analysis and uncertainty propagation for robust design, SAE transactions
  (2004) 121--128.

\bibitem{Horiguchi21}
A.~Horiguchi, M.~T. Pratola, T.~J. Santner, {Assessing variable activity for
  Bayesian regression trees}, Reliability Engineering \& System Safety 207
  (2021) 107391.

\bibitem{Antoniadis21}
A.~Antoniadis, S.~Lambert-Lacroix, J.-M. Poggi, {Random forests for global
  sensitivity analysis: A selective review}, Reliability Engineering \& System
  Safety 206 (2021) 107312.

\bibitem{Steiner19}
M.~Steiner, J.-M. Bourinet, T.~Lahmer, An adaptive sampling method for global
  sensitivity analysis based on least-squares support vector regression,
  Reliability Engineering \& System Safety 183 (2019) 323--340.

\bibitem{Fock14}
E.~Fock, {Global Sensitivity Analysis Approach for Input Selection and System
  Identification Purposes—A New Framework for Feedforward Neural Networks},
  IEEE Transactions on Neural Networks and Learning Systems 25~(8) (2014)
  1484--1495.

\bibitem{Datteo18}
A.~Datteo, G.~Busca, G.~Quattromani, A.~Cigada, {On the use of AR models for
  SHM: A global sensitivity and uncertainty analysis framework}, Reliability
  Engineering \& System Safety 170 (2018) 99--115.

\bibitem{Leak14}
E.~Todri, A.~Amenaghawon, I.~{Del Val}, D.~Leak, C.~Kontoravdi, S.~Kucherenko,
  N.~Shah, Global sensitivity analysis and meta-modeling of an ethanol
  production process, Chemical Engineering Science 114 (2014) 114--127.

\bibitem{Kai20}
K.~Cheng, L.~Zhenzhou, C.~Ling, S.~Zhou, Surrogate-assisted global sensitivity
  analysis: an overview, Structural and Multidisciplinary Optimization 61
  (2020).

\bibitem{Wu16}
Z.~Wu, D.~Wang, P.~N. Okolo, F.~Hu, W.~Zhang, {Global sensitivity analysis
  using a Gaussian Radial Basis Function metamodel}, Reliability Engineering \&
  System Safety 154 (2016) 171--179.

\bibitem{Wu19}
Z.~Wu, W.~Wang, D.~Wang, K.~Zhao, W.~Zhang, Global sensitivity analysis using
  orthogonal augmented radial basis function, Reliability Engineering \& System
  Safety 185 (2019) 291--302.

\bibitem{Blatman10}
G.~Blatman, B.~Sudret, Efficient computation of global sensitivity indices
  using sparse polynomial chaos expansions, Reliability Engineering \& System
  Safety 95~(11) (2010) 1216--1229.

\bibitem{Luthen21}
N.~L\"{u}then, S.~Marelli, B.~Sudret, {Sparse Polynomial Chaos Expansions:
  Literature Survey and Benchmark}, SIAM/ASA Journal on Uncertainty
  Quantification 9~(2) (2021) 593--649.

\bibitem{AlexanderianGremaudSmith20}
A.~Alexanderian, P.~A. Gremaud, R.~C. Smith, Variance-based sensitivity
  analysis for time-dependent processes, Reliability Engineering \& System
  Safety 196 (2020) 106722.

\bibitem{Ehre20}
M.~Ehre, I.~Papaioannou, D.~Straub, {Global sensitivity analysis in high
  dimensions with PLS-PCE}, Reliability Engineering \& System Safety 198 (2020)
  106861.

\bibitem{Zhou20}
Y.~Zhou, Z.~Lu, J.~Hu, Y.~Hu, Surrogate modeling of high-dimensional problems
  via data-driven polynomial chaos expansions and sparse partial least square,
  Computer Methods in Applied Mechanics and Engineering 364 (2020) 112906.

\bibitem{Luthen22}
N.~L\"{u}then, S.~Marelli, B.~Sudret, Automatic selection of basis-adaptive
  sparse polynomial chaos expansions for engineering applications,
  International Journal for Uncertainty Quantification 12~(3) (2022) 49--74.

\bibitem{Almohammadi22}
S.~M. Almohammadi, O.~P.~L. Maître, O.~M. Knio, {Computational Challenges in
  Sampling and Representation of Uncertain Reaction Kinetics in Large
  Dimensions}, International Journal for Uncertainty Quantification 12~(1)
  (2022) 1--24.

\bibitem{Zhang22}
Q.~Zhang, Y.-G. Zhao, K.~Kolozvari, L.~Xu, Reliability analysis of reinforced
  concrete structure against progressive collapse, Reliability Engineering \&
  System Safety 228 (2022) 108831.

\bibitem{HuangELM}
G.-B. Huang, Q.-Y. Zhu, C.-K. Siew, {Extreme learning machine: Theory and
  applications}, Neurocomputing 70~(1) (2006) 489--501, neural Networks.

\bibitem{Huang11}
G.-B. Huang, D.~Wang, Y.~Lan, Extreme learning machines: a survey,
  International Journal of Machine Learning and Cybernetics 2~(2) (2011)
  107--122.

\bibitem{Stewart21}
P.~Jorgensen, D.~E. Stewart, {Approximation Properties of Ridge Functions and
  Extreme Learning Machines}, SIAM Journal on Mathematics of Data Science 3~(3)
  (2021) 815--832.

\bibitem{Huang15}
G.~Huang, G.-B. Huang, S.~Song, K.~You, {Trends in extreme learning machines: A
  review}, Neural Networks 61 (2015) 32--48.

\bibitem{Wang22}
J.~Wang, S.~Lu, S.-H. Wang, Y.-D. Zhang, A review on extreme learning machine,
  Multimedia Tools and Applications 81~(29) (2022) 41611--41660.

\bibitem{Nagawkar21}
J.~Nagawkar, L.~Leifsson, {Efficient Global Sensitivity Analysis of Model-Based
  Ultrasonic Nondestructive Testing Systems Using Machine Learning and Sobol’
  Indices}, Journal of Nondestructive Evaluation, Diagnostics and Prognostics
  of Engineering Systems 4~(4) (2021).

\bibitem{Walzberg21}
J.~Walzberg, A.~Carpenter, G.~A. Heath, Role of the social factors in success
  of solar photovoltaic reuse and recycle programmes, Nature Energy 6~(9)
  (2021) 913--924.

\bibitem{Li19}
S.~Li, B.~Yang, F.~Qi, {Accelerate global sensitivity analysis using artificial
  neural network algorithm: Case studies for combustion kinetic model},
  Combustion and Flame 168 (2016) 53--64.

\bibitem{Kapusuzoglu21}
B.~Kapusuzoglu, S.~Mahadevan, Information fusion and machine learning for
  sensitivity analysis using physics knowledge and experimental data,
  Reliability Engineering \& System Safety 214 (2021) 107712.

\bibitem{Ye21}
D.~Ye, A.~Nikishova, L.~Veen, P.~Zun, A.~G. Hoekstra, Non-intrusive and
  semi-intrusive uncertainty quantification of a multiscale in-stent restenosis
  model, Reliability Engineering \& System Safety 214 (2021) 107734.

\bibitem{Zhao23}
Y.~Zhao, X.~Cheng, T.~Zhang, L.~Wang, W.~Shao, J.~Wiart, A global–local
  attention network for uncertainty analysis of ground penetrating radar
  modeling, Reliability Engineering \& System Safety 234 (2023) 109176.

\bibitem{Engelbrecht01}
A.~Engelbrecht, A new pruning heuristic based on variance analysis of
  sensitivity information, Neural Networks, IEEE Transactions on 12 (2001) 1386
  -- 1399.

\bibitem{Blalock2020}
D.~Blalock, J.~J. Gonzalez~Ortiz, J.~Frankle, J.~Guttag, What is the state of
  neural network pruning?, Proceedings of machine learning and systems 2 (2020)
  129--146.

\bibitem{Pinkus99}
A.~Pinkus, {Approximation theory of the MLP model in neural networks}, Acta
  Numerica 8 (1999) 143–195.

\bibitem{Leshno93}
M.~Leshno, V.~Y. Lin, A.~Pinkus, S.~Schocken, Multilayer feedforward networks
  with a nonpolynomial activation function can approximate any function, Neural
  Networks 6~(6) (1993) 861--867.

\bibitem{Hansen10}
P.~C. Hansen, {Getting Serious: Choosing the Regularization Parameter}, Society
  for Industrial and Applied Mathematics, 2010, Ch.~5, pp. 85--107.

\bibitem{Saltelli10}
A.~Saltelli, P.~Annoni, I.~Azzini, F.~Campolongo, M.~Ratto, S.~Tarantola,
  {Variance based sensitivity analysis of model output. Design and estimator
  for the total sensitivity index}, Computer Physics Communications 181~(2)
  (2010) 259--270.

\bibitem{Vilar02}
J.~M. Vilar, H.~Y. Kueh, N.~Barkai, S.~Leibler, {Mechanisms of noise-resistance
  in genetic oscillators}, Proceedings of the National Academy of Sciences
  99~(9) (2002) 5988--5992.

\bibitem{Sheppard12}
P.~W. Sheppard, M.~Rathinam, M.~Khammash, A pathwise derivative approach to the
  computation of parameter sensitivities in discrete stochastic chemical
  systems, The Journal of Chemical Physics 136~(3) (2012) 034115.

\bibitem{Merritt21}
M.~Merritt, A.~Alexanderian, P.~A. Gremaud, {Multiscale Global Sensitivity
  Analysis for Stochastic Chemical Systems}, Multiscale Modeling \& Simulation
  19~(1) (2021) 440--459.

\bibitem{Owen13}
A.~B. Owen, \href{https://statweb.stanford.edu/~owen/mc/}{{Monte Carlo theory,
  methods and examples}}, Art B. Owen, 2013.
\newline\urlprefix\url{https://statweb.stanford.edu/~owen/mc/}

\end{thebibliography}

\appendix
\section{Derivation of analytic formulas for ELM surrogate}
\label{sec:derivation}
\subsection{Proof of \cref{lma1}}
\begin{proof} For the surrogate presented in this paper, we will derive expressions for the mean and variance. 
Recall the expression for the ELM \eqref{equ:elm} and 

the definition of the function $\epsilon(t)$, 
given in the statement of the lemma,

\[
\epsilon(t) := \int_0^1e^{tx}dx 
= \left\{\begin{array}{ccc}
\frac{e^t-1}{t} && t\neq 0, \\
1 && t = 0.
\end{array}\right. 
\]

%~~~~~~~~~~~~~~~~~~ Mean

We first find the expression for the mean: 
\[
\E(\hat{f})=\int_{[0,1]^d}\sum_{j=1}^n\Big(\beta_je^{b_j}\prod_{l=1}^de^{w_{j,l}x_l}\Big)d\vec{x}=\sum_{j=1}^n\Big(\beta_je^{b_j}\prod_{l=1}^d\int_0^1e^{w_{j,l}x_l}dx_l\Big).
\]
%Note the definition of the function $\epsilon(t)$. 
Therefore, the mean can be expressed as 
\[\E(\hat{f})=\sum_{j=1}^n\Big(\beta_je^{b_j}\prod_{l=1}^d\epsilon(w_{j,l})\Big).
\]

%%~~~~~~~~~~~~~~ Variance

Next, we find the expression for the variance: 
\[
\begin{aligned}
\var(\hat{f})
=
\E(\hat{f}^2)-\E(\hat{f})^2
&=
\int_{[0,1]^d}\Big( \sum_{j=1}^n\big(\beta_je^{b_j}\prod_{l=1}^de^{w_{j,l}x_l}\big)\Big)^2d\vec{x}-\E(\hat{f})^2
\\
&=
\sum_{j,i=1}^n\Big(\beta_j\beta_ie^{b_j+b_i}\prod_{l=1}^d\epsilon(w_{j,l}+w_{i,l})\Big)-\E(\hat{f})^2\\
&=
\sum_{j,i=1}^n\beta_j\beta_ie^{b_j+b_i}\Big(\prod_{l=1}^d\epsilon(w_{j,l}+w_{i,l})-\prod_{r=1}^d\epsilon(w_{j,r})\epsilon(w_{i,r})\Big).
\end{aligned}
\]
\end{proof}

\subsection{Proof of \cref{prop1}}
\begin{proof}
For the surrogate presented in this paper, we will derive expressions for the general formulas for regular and total Sobol' indices and offer simplified expressions for indices corresponding to single variables. 
Recall the expression for the ELM~\eqref{equ:elm} and the function $\epsilon(t)$, and 
let $\E(\hat{f})$ and $\var(\hat{f})$ be as presented in \cref{lma1}. 

Given subset of variables of $\{x_1,\dots,x_d\}$, we first derive the expressions for Sobol' indices,
\[
S_\vec{u}=\frac{\var(\f_\vec{u})}{\var(\f)},\quad \f_\vec{u}:=\sum_{\vec{v}\subseteq\vec{u}}(-1)^{|\vec{u}|-|\vec{v}|}\E(\f|x_l,\ l\in\vec{v}), 
\]
where $\vec{u}\subset\{1,\dots,d\}$~\cite{PrieurTarantola17}. 
Since terms in the ANOVA decomposition have the property that $\var(\hat{f}_\vec{u})=\E(\hat{f}_\vec{u}^2)$~\cite{Owen13} we can we can express the Sobol' index $S_\vec{u}$ for~\eqref{equ:elm} as
\[
\begin{aligned}
S_\vec{u} &= \frac{\var(\hat{f}_\vec{u})}{\var(\hat{f})}=\frac{\E(\hat{f}_\vec{u}^2)}{\var(\hat{f})}\\ 
&=
\frac{1}{\var(\hat{f})}\int_{[0,1]^d}\Big(\sum_{\vec{v}\subseteq\vec{u}}\sum_{j=1}^n(-1)^{|\vec{u}|-|\vec{v}|}\beta_je^{b_j}\Big(\prod_{l\in\vec{v}}e^{w_{j,l}x_l}\Big)\Big(\prod_{r\not\in\vec{v}}\epsilon(w_{j,r})\Big)\Big)^2d\vec{x}
\\
&=
\frac{1}{\var(\hat{f})}\sum_{\vec{v},\tilde{\vec{v}}\subseteq\vec{u}}\sum_{j,i=1}^n(-1)^{2|\vec{u}|-|\vec{v}|-|\tilde{\vec{v}}|}\int_{[0,1]^d}\Big(\beta_j\beta_ie^{b_j+b_i}\prod_{l\in\vec{v}}e^{w_{j,l}x_l}\prod_{r\not\in\vec{v}}\epsilon(w_{j,r})\prod_{s\in\tilde{\vec{v}}}e^{w_{i,s}x_l}\prod_{q\not\in\tilde{\vec{v}}}\epsilon(w_{i,q})\Big)d\vec{x} 
\\
&=
 \frac{1}{\var(\hat{f})}\sum_{\vec{v},\tilde{\vec{v}}\subseteq\vec{u}}\sum_{j,i=1}^n(-1)^{2|\vec{u}|-|\vec{v}|-|\tilde{\vec{v}}|}\beta_j\beta_ie^{b_j+b_i}\Big(\prod_{l\in\vec{v}\cap\tilde{\vec{v}}}\epsilon(w_{j,l}+w_{i,l})\Big)\Big(\prod_{r\not\in\vec{v}\cap\tilde{\vec{v}}}\epsilon(w_{j,r})\epsilon(w_{i,r})\Big).
\end{aligned}
\]
When we only consider the single variable case or, in other words, when $\vec{u}=\{k\}$, then we arrive at equation~\eqref{equ:reg} for the first-order Sobol' indices:
\[
\begin{aligned}
S_k &= \frac{1}{\var(\f)}\sum_{j,i=1}^n\Big(\beta_j\beta_ie^{b_j+b_i}\epsilon(w_{j,k}+w_{i,k})\prod_{l\neq k}\epsilon(w_{j,l})\epsilon(w_{i,l})\Big)-\frac{\E(\f)^2}{\var(\f)}\\
&=\frac{1}{\var(\hat{f})}\sum_{j,i=1}^n\beta_j\beta_ie^{b_j+b_i}
\big(\epsilon(w_{j,k}+w_{i,k})-\epsilon(w_{j,k})\epsilon(w_{i,k})\big)\prod_{l\neq k}\epsilon(w_{j,l})\epsilon(w_{i,l}).
\end{aligned}
\]

 Now we derive the expressions for the total Sobol' indices for a given subset of variables of $\{x_1,\dots,x_d\}$,
\[
S^\mathrm{tot}_\vec{u}=1-\frac{\var(\E(\f|x_r,\  r\not\in\vec{u}))}{\var(\f)}, 
\]
where $\vec{u}\subset\{1,\dots,d\}$~\cite{PrieurTarantola17}. We find the expression for the total index $S^\mathrm{tot}_\vec{u}$ for~\eqref{equ:elm}:
\[
\begin{aligned}
S^\mathrm{tot}_\vec{u} &= 1 - \frac{1}{\var(\f)}\Big(\int_{[0,1]^d}\Big(\sum_{j=1}^n\beta_je^{b_j}\Big(\prod_{l\in\vec{u}}\epsilon(w_{j,l})\Big)\Big(\prod_{r\not\in\vec{u}}e^{w_{j,r}x_r})\Big)\Big)^2d\vec{x}-\E(\f)^2\Big)
\\
&=
1 - \frac{1}{\var(\f)}\Big(\sum_{j,i=1}\beta_j\beta_ie^{b_j+b_i}\Big(\prod_{l\in\vec{u}}\epsilon(w_{j,l})\epsilon(w_{i,l})\Big)\Big(\prod_{r\not\in\vec{u}}\epsilon(w_{j,r}+w_{i,r})\Big)-\E(\f)^2\Big).
\end{aligned}
\]
When we consider the total Sobol' index corresponding to a single variable $x_k$, we have expression~\eqref{equ:tot}:
\[
\begin{aligned}
S_k^{\mathrm{tot}} &= 1- \frac{1}{\var(\hat{f})}\Big(\sum_{j,i=1}^n\beta_j\beta_ie^{b_j+b_i}\epsilon(w_{j,k})\epsilon(w_{i,k})\prod_{l\neq k}\epsilon(w_{j,l}+w_{i,l})-\E(\hat{f})^2\Big) 
\\
&= 
1- \frac{1}{\var(\hat{f})}\sum_{j,i=1}^n\beta_j\beta_ie^{b_j+b_i}\epsilon(w_{j,k})\epsilon(w_{i,k})\Big(\prod_{l\neq k}\epsilon(w_{j,l}+w_{i,l})-\prod_{r\neq k}\epsilon(w_{j,r})\epsilon(w_{i,r})\Big).
\end{aligned}
\]
\end{proof}

\section{Sobol' indices of analytic example~\eqref{equ:bnch_int}}\label{sec:sobolinds}

Here we present formulas for the analytic first-order and total Sobol' indices for~\eqref{equ:bnch_int}, given an input dimension $d$ and parameter $\delta$. 
The first-order Sobol indices are the same for all $k=1,\dots,d$. The total
Sobol' indices are also the same for $k=1,\dots,d$. The first-order and total
Sobol' indices are given by
\[
S_k=\frac{\var(\E( f_\delta | x_k))}{\var(f_\delta)}\quad\text{and}\quad S_k^\mathrm{tot}=1-\frac{\var(\E(f_\delta|x_l,\ l\neq k))}{\var(f_\delta)},
\]
where
\[
\begin{aligned}
&\var(f_\delta)=\frac{d\delta}{9}\Big(\frac{3}{2}\Big)^d+\delta^2\Big(\Big(\frac{7}{3}\Big)^d-\Big(\frac{9}{4}\Big)^d\Big)+\frac{d}{12}, \\ 
&\var(\E( f_\delta | x_k))=\frac{\delta^2}{27}\Big(\frac{9}{4}\Big)^d+\frac{\delta}{9}\Big(\frac{3}{2}\Big)^d+\frac{1}{12},\\ 
&\var(\E(f_\delta|x_l,\ l\neq k))=\frac{\delta(d-1)}{9}\Big(\frac{3}{2}\Big)^d+\delta^2\Big(\frac{27}{28}\Big(\frac{7}{3}\Big)^d-\Big(\frac{9}{4}\Big)^d\Big)+\frac{d-1}{12}. 
\end{aligned}
\]

\section{Setup of the model in~\cref{sec:highdim}}
\label{sec:highdim_details}
Here we describe our choices for $\vec{x}_0$ and $\mat{Q}$
in~\eqref{equ:eigsol}, as implemented in \MATLAB.  For the initial state, $\boldsymbol x_0$, we
take $\boldsymbol x_0=[1\quad\dots\quad 1]^\top$.  We chose $\mat Q$ by
generating a random orthogonal matrix as follows.  First, we generated a matrix
$\mat{P}$ whose entries are independent draws from standard normal
distribution. This was done using \MATLAB's \verb+randn+ command. The random
seed was fixed using \verb+rag(1)+.  Subsequently, the matrix $\mat Q$ was
obtained by performing a QR factorization $\mat{P} = \mat{Q} \mat{R}$, where
$\mat{Q}$ is orthogonal and $\mat{R}$ is upper triangular, and retaining
$\mat{Q}$.  The QR factorization was computed using \MATLAB's \verb+qr+
command.

\end{document}